\documentclass[11pt]{ijnam}
\def\currentvolume{1}
\def\currentissue{1}
\def\currentyear{2004}
\def\month{March}
\pagespan{1}{18}
\copyrightinfo{2004}{} 

\usepackage{stmaryrd}
\newcommand{\jump}[1]{\left\llbracket#1\right\rrbracket} 

\newcommand{\commentout}[1]{{}} 

\renewcommand{\theequation}{\arabic{section}.\arabic{equation}}

\newtheorem{lemma}{Lemma}[section]

\newcommand{\normtvb}[1]{{\left\vert\kern-0.25ex\left\vert\kern-0.25ex\left\vert #1
    \right\vert\kern-0.25ex\right\vert\kern-0.25ex\right\vert}}

\newcommand{\cN}{\mathcal{N}}

\newcommand{\bfn}{{\bf n}}

\newcommand{\hK}{{\hat {K}}}

\newcommand{\hu}{{\hat {u}}}

\newcommand{\hv}{{\hat {v}}}

\newcommand{\hbeta}{{\hat{\beta}}}

\newcommand{\hGamma}{{\hat{\Gamma}}}

\newcommand{\hphi}{{\hat{\phi}}}

\usepackage{bm}
\usepackage{graphicx}
\usepackage{mathrsfs}
\usepackage{booktabs}
\usepackage{multirow}
\usepackage{color}
\theoremstyle{remark}\newtheorem{remark}{Remark}
\numberwithin{equation}{section}

\newcommand{\xz}[1]{{\color{black} #1}}

\begin{document}

\title[Frenet IFE on Triangular Meshes]
{Frenet Immersed Finite Element Spaces on Triangular Meshes}

\author{Yuanhui Lin}
\address{
  School of Mathematics, Sichuan University, Chengdu, Sichuan 610065, China.
}
\email{l1006569522@outlook.com}

\author{Xu Zhang}
\address{
  Department of Mathematics, Oklahoma State University, Stillwater, OK 74078, USA.
}
\email{xzhang@okstate.edu}

\author{Tao Lin\quad}
\address{Department of Mathematics, Virginia Tech, Blacksburg, VA 24060, USA.}
\email{tlin@vt.edu}


\date{December 1, 2025 and, April 25, 2026 in revised form.}

\thanks{}


\abstract{In this paper, we develop geometry-conforming immersed finite element (IFE) spaces on triangular meshes for elliptic interface problems. The construction is built on a Frenet–Serret mapping that transforms a smooth interface curve into a straight line, \xz{so that the interface jump conditions can be imposed exactly}. Extending the framework of \cite{2025AdjeridLinMeghaichi} from rectangular meshes to triangular meshes, we introduce three types of high-order Frenet-IFE constructions: an initial construction using monomial bases, a \xz{general} construction using orthogonal polynomials, and reconstructed IFE \xz{bases} designed to improve the conditioning of the mass matrix. The approximation properties of these new IFE spaces are investigated through extensive numerical experiments. We also incorporate the new IFE spaces into interior penalty discontinuous Galerkin methods for solving elliptic interface problems, and demonstrate optimal convergence rates in $H^1$- and $L^2$- norms.}

\maketitle

\section{Introduction}
In this paper, we develop a high-order geometry-conforming immersed finite element (GC-IFE) on triangular meshes for the second-order elliptic interface problem: 
\begin{eqnarray}
-{\nabla}\cdot(\beta\nabla u) = {f}, &\mbox{in}& \Omega^- \cup \Omega^+, \label{eq:FrenentIFE_ellipPDE_tri}\\
{u} = {g}, &\mbox{on}& \partial \Omega. \label{eq:FrenentIFE_ellipPDE_BC_tri} 
\end{eqnarray}
Here, $\Omega\subset \mathbb{R}^2$ is a polygonal domain separated by a smooth interface $\Gamma$ into two subdomains $\Omega^-$ and $\Omega^+$. The coefficient function $\beta$ is discontinuous across the interface. Without loss of generality, we assume that it is a piecewise constant $\beta|_{\Omega^\pm} = \beta^\pm$. The solution $u$ is assumed to satisfy the following interface jump conditions: 
\begin{eqnarray}
\jump{u}_\Gamma = 0, &\mbox{on}~\Gamma, \label{eq: jump1} \\
\jump{\beta\frac{\partial u}{\partial \bfn}}_\Gamma = 0, &\mbox{on} ~\Gamma, \label{eq: jump2}
\end{eqnarray}
where $\jump{v}_\Gamma :=v^+|_\Gamma-v^-|_\Gamma$ is the jump across the interface, and $\mathbf{n}$ denotes the unit normal on $\Gamma$ from $\Omega^-$ to $\Omega^+$. As described in \cite{2017AdjeridGuoLin, 2024AdjeridLinMeghaichi, 2019GuoLin2}, with the intention of constructing a high-order $\mathbb{P}_m$ IFE space, we further assume that the solution $u$ satisfies the Laplacian extended jump conditions:
\begin{eqnarray}
\jump{ \beta\frac{\partial^{j} \triangle u}{\partial\mathbf{ n}^{j}} }_{\Gamma} = 0, \quad j= 0, 1, \cdots, m-2.
\label{eq: jumpEx}
\end{eqnarray}

The interface problem \eqref{eq:FrenentIFE_ellipPDE_tri} - \eqref{eq: jumpEx} arises widely in  science and engineering applications. The Immersed Finite Element (IFE) method \cite{2019GuoLinZhuang,2021GuoZhang,2026JiLi,2023JiWangChenLi,2005KafafyLinLinWang,1998Li,2015LinLinZhang,2019LinSheenZhang} is a class of numerical techniques designed to solve such problems without requiring the mesh to align with the interface. Its key idea  is to modify the local shape functions only on interface elements according to the prescribed jump conditions, while leaving the finite element space unchanged on regular elements. 

Traditionally, local IFE functions are constructed as piecewise polynomials on subelements separated by the interface. For low-order IFE spaces \cite{2019GuoLin,2008HeLinLin,2022JiWangChenLi2,2004LiLinLinRogers}, the interface can be adequately represented by a linear approximation, and the associated IFE spaces are built accordingly. 
However, for high-order IFE spaces, such a linear approximation becomes insufficient. When the interface is an arbitrary smooth curve, the jump conditions cannot be enforced exactly because polynomial pieces cannot match perfectly along the interface curve. As a result, the jump conditions are often imposed on the interface curve only in a weak sense. Examples include least-squares enforcement of the jump conditions \cite{2017AdjeridGuoLin, 2021ChenZhang, 2024ChenZhang, 2025ChenZhang},  Cauchy extension \cite{2023AdjeridBabuskaGuoLin, 2019GuoLin2}, and \xz{enforcement of jump conditions at selected points} \cite{2022Adjerid,2014AdjeridBenromdhaneLin, 2017GuzmanSanchezSarkis}. While these approaches yield IFE spaces with sufficient approximation capability, the resulting spaces are geometrically nonconforming, meaning the local IFE functions generally do not lie in $H^1$. This drawback has two main ramifications. First, an appropriate penalty must be incorporated along the interface to mitigate the effect of the inherent discontinuity of the IFE functions \cite{2023AdjeridBabuskaGuoLin, 2026ZhuangZhangSarkisLin}, which consequently increases the complexity of the IFE scheme. Second, the coexistence of discontinuity and penalty terms along the interface poses additional challenges for the associated error analysis. 

Recently, a geometrically conforming IFE space was introduced on rectangular meshes \cite{2024AdjeridLinMeghaichi}. This construction uses the Frenet-Serret transformation, an idea from differential geometry, to map an arbitrary interface curve in the Cartesian $x$-$y$ plane into a straight line segment in the Frenet coordinates $\eta$-$\xi$. Although standard tensor-product polynomial spaces $\mathbb{Q}_m$ are used on Frenet reference elements, the corresponding IFE shape functions on the physical element are no longer piecewise polynomials because of the nonlinearity of the Frenet mapping. Nevertheless, these new IFE functions preserve continuity inside each interface element, which makes the resulting space geometry-conforming. In \cite{2025AdjeridLinMeghaichi}, the approximation properties of the IFE spaces were theoretically established and the optimal {\it a priori} error estimates for immersed discontinuous Galerkin solution were proved. In \cite{2025AdjeridLinMeghaichi2}, some reconstruction techniques for Frenet-IFE bases were further developed to improve the conditioning of the associated mass matrix.

The Frenet-IFE framework \cite{2024AdjeridLinMeghaichi,2025AdjeridLinMeghaichi2,2025AdjeridLinMeghaichi} has so far been developed exclusively for structured domains discretized by rectangular meshes. This restriction limits its applicability, since many practical problems involve domains with complex geometries. In particular, rectangular meshes cannot be constructed on non-rectangular polygonal domains or domains with curved boundaries. In contrast, triangular meshes provide much greater flexibility. They can naturally conform to complex boundaries, and enable efficient local mesh refinement. These advantages motivate us to generalize and extend the Frenet-IFE framework from Cartesian rectangular meshes to unstructured triangular meshes. 


This paper is strongly influenced by, and builds upon, the developments in \cite{2025AdjeridLinMeghaichi2}. We present three approaches for constructing high-order, geometrically conforming IFE spaces on triangular meshes. The first approach, referred to as  {\it initial construction}, uses $(m+1)(m+2)/2$ monomials forming the basis of $\mathbb{P}_m$. Thanks to its simple structure, $m(m+1)/2$ of the resulting IFE basis functions can be written down explicitly, while the remaining $m+1$ functions are determined by solving a small linear system. The second approach, called {\it general construction}, uses orthogonal polynomials to build IFE basis functions, providing a more flexible framework. The third approach, named {\it reconstruction}, modifies the Frenet IFE basis to achieve optimally conditioned mass matrices, which is crucial for numerical stability, especially for time dependent problems. 

The remainder of the paper is organized as follows. Section \ref{sec: prelim} reviews the necessary notations, the Frenet transformation, and related results previously established for rectangular meshes. Sections \ref{sec: initial construction} and \ref{sec: general construction basis} present the initial and general constructions of the Frenet-IFE basis functions on triangular meshes. In Section \ref{sec: reconstruction basis}, we apply the reconstruction strategies \cite{2025AdjeridLinMeghaichi2} to IFE basis functions on triangular meshes for improving the conditioning of mass matrices.
In Section \ref{sec: numerical}, we report numerical results for the $L^2$ projection and discontinuous Galerkin approximation of the elliptic interface problem. Finally, Section \ref{sec: conclusion} provides some concluding remarks.

\section{Preliminaries}\label{sec: prelim}
In this section, we review some notation for Frenet transformation and recall related results for the rectangular meshes previously presented in \cite{2024AdjeridLinMeghaichi,2025AdjeridLinMeghaichi2,2025AdjeridLinMeghaichi}. 

Let $\mathcal{T}_h = \{K\}$ be a shape-regular triangulation of $\Omega$ that is independent of the interface $\Gamma$. We call $K$ an interface element if $K\cap \Gamma \xz{\neq} \emptyset$, and denote by $\mathcal{T}_h^i$ the collection of all interface elements. The set of noninterface elements is then $\mathcal{T}_h^n = \mathcal{T}_h \setminus\mathcal{T}_h^i$.

Let $\Gamma$ be a smooth interface in $\mathbb{R}^2$, parametrized by
\begin{equation}\label{eq: g}
\mathbf{g}(\xi) = [g_1(\xi),g_2(\xi)]: [\xi_s,\xi_e]\to \mathbb{R}^2.
\end{equation}
The Frenet–Serret apparatus associated with $\mathbf{g}(\xi)$ consists of the unit tangent vector $\bm{\tau}(\xi)$, the unit normal vector 
$\bm{n}(\xi)$, and the curvature $\kappa(\xi)$ defined by 
\begin{equation}\label{eq: Frenet apparatus}
\bm{\tau}(\xi) = \frac{\mathbf{g}'(\xi)}{\|\mathbf{g}'(\xi)\|},\quad
\bm{n}(\xi)  = Q\bm{\tau}(\xi),\quad
\kappa(\xi)  =  \frac{\mathbf{g}'(\xi)^TQ\mathbf{g}''(\xi)}{\|\mathbf{g}'(\xi)\|^3}
\end{equation}
where
\[Q = 
\begin{bmatrix}
0 & 1 \\-1 & 0
\end{bmatrix}.\]
The mapping $P_\Gamma: (\eta,\xi) \to (x,y)$ is defined by 
\[
\mathbf{x}(\eta,\xi) = \begin{bmatrix}
x(\eta,\xi) \\y(\eta,\xi)
\end{bmatrix} = P_\Gamma(\eta,\xi) = \mathbf{g}(\xi)+\eta\bm{n}(\xi).
\]
It is well-known \cite{2012AbateTovena} that $\Gamma$ has an $\epsilon$-tubular neighborhood 
\[N_\Gamma(\epsilon) = P_\Gamma([-\epsilon,\epsilon]\times[\xi_s,\xi_e])\] 
within which the transform $P_\Gamma$ is a bijection. Hence its inverse 
\begin{equation}
    R_\Gamma = P_\Gamma^{-1}: N_\Gamma(\epsilon) \to 
    [-\epsilon,\epsilon]\times[\xi_s,\xi_e]
\end{equation} 
is well defined 
such that 
\[
\begin{bmatrix} \eta \\ \xi \end{bmatrix}  
= 
\begin{bmatrix} \eta(x,y) \\ \xi(x,y) \end{bmatrix}  
 = R_\Gamma(x,y).
 \]
For sufficiently small mesh size $h$, every interface element in $\mathcal{T}_h$ lies inside $N_\Gamma(\epsilon)$. Consider an interface element $K\in\mathcal{T}_h^i$ with vertices $ A_1$, $A_2$, and $A_3$, and let $R_\Gamma(A_i)= [\eta_i,\xi_i] $. \xz{Let $\xi_{0,K} = \min\limits_{1\le i\le 3}\{\xi_i\}$ and  $\xi_{1,K} = \max\limits_{1\le i\le 3}\{\xi_i\}$.}  We form a fictitious element $K_F$ containing $K$, with two curved edges parallel to the interface $\Gamma$ and two straight edges. By applying the inverse mapping, the interface triangle $K$ becomes a curved triangle $\hat K = R_\Gamma(K) \subset \widehat{K}_F$. A key feature of the Frenet mapping is that the interface segment $\Gamma_{K_F} = \Gamma\cap K_F$ is mapped to the vertical line segment $\widehat{\Gamma}_{K_F} = R_\Gamma(\Gamma_{K_F})$ in the $(\eta,\xi)$-plane. 

\begin{figure}[h]
\includegraphics[width = .7\textwidth]{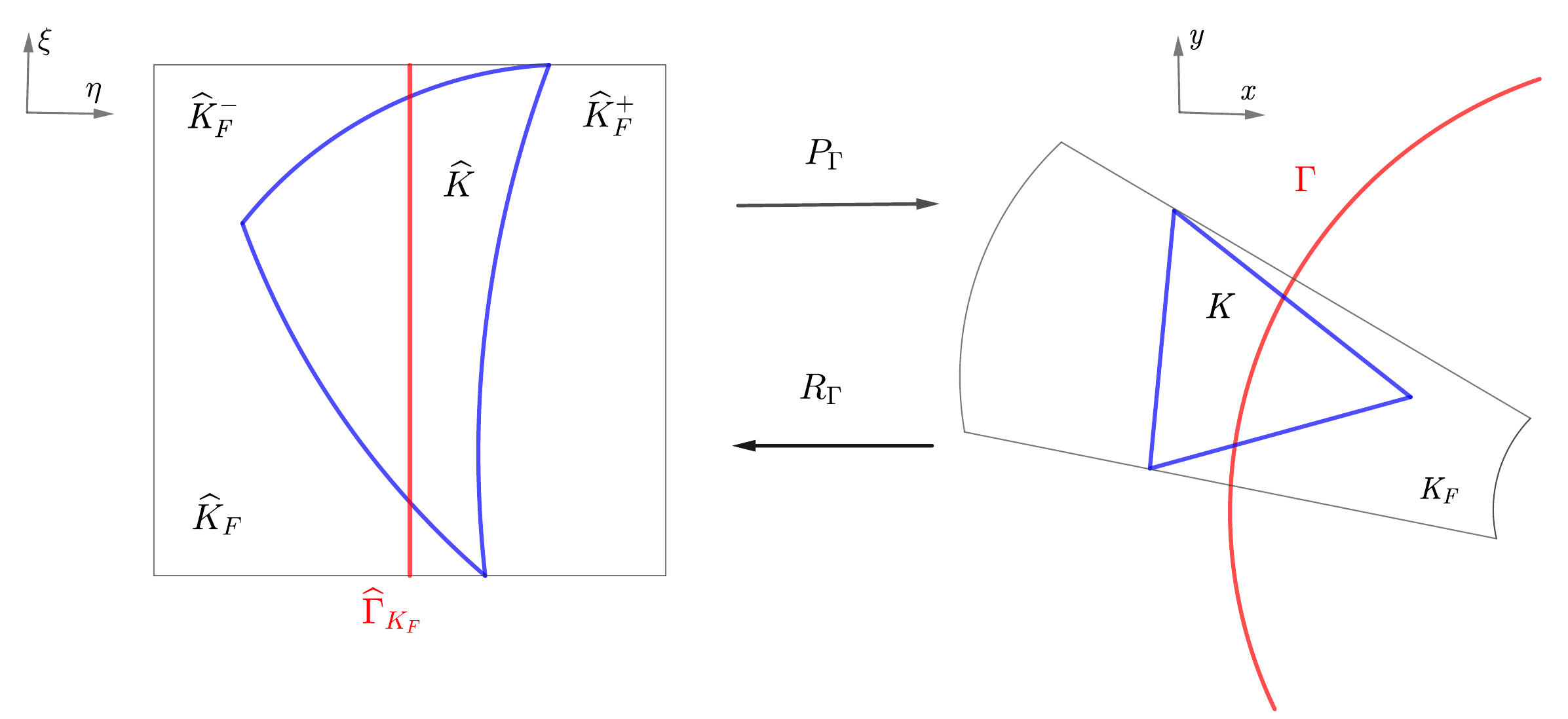}
\caption{An illustration of an interface triangle $K$ and its Frenet-mapping.}
\label{fig: mapping}
\end{figure}

See Figure \ref{fig: mapping} for an illustration. 

Using this coordinate transformation, every function $v: N_\Gamma(\epsilon) \rightarrow \mathbb{R}$ is associated  with the following function of the Frenet coordinates:
\begin{eqnarray}
\hv = v\circ P_\Gamma: [-\epsilon, \epsilon]\times [\xi_s, \xi_e] \to \mathbb{R}. \label{eq:uhat}
\end{eqnarray}
Moreover, the interface jump conditions in \eqref{eq: jump1}-\eqref{eq: jumpEx} are transformed to
\begin{eqnarray}
\jump{\hat u}_{\widehat{\Gamma}_{K_F}} &=& 0,\qquad \label{eq: Fr jump1}\\
\jump{\hat\beta\hat u_\eta}_{\widehat{\Gamma}_{K_F}} &=&0,\qquad\label{eq: Fr jump2}\\
\jump{\hat\beta\partial_\eta^j \mathcal{L}(\hat u)}_{\widehat{\Gamma}_{K_F}}  &=&0,\quad j = 0,1,\cdots,m-2,\label{eq: Fr jump3}
\end{eqnarray}
where the Laplacian in Frenet coordinate is given by
\begin{equation}\label{eq: lap}
\mathcal{L}(\hat u(\eta,\xi)) = \hat u_{\eta\eta} + J_0(\eta,\xi)\hat u_{\xi\xi}(\eta,\xi)+J_1(\eta,\xi)\hat u_{\eta}(\eta,\xi) + J_2(\eta,\xi)\hat u_{\xi}(\eta,\xi)
\end{equation}
with
\begin{equation}\label{eq: J}
\begin{split}
J_0(\eta,\xi) &= \rho^2(\eta,\xi),\quad
J_1(\eta,\xi) = \kappa(\xi)\psi(\eta,\xi), \\
J_2(\eta,\xi) &= \rho^2(\eta,\xi)\left(\eta \kappa'(\xi)\psi(\eta,\xi)  \frac{\mathbf{g}'(\xi)\cdot \mathbf{g}''(\xi)}{\|\mathbf{g}'(\xi)\|^2}\right).
\end{split}
\end{equation}
\xz{
and 
\[
\rho(\eta,\xi) = \|\mathbf{g}'(\xi)\|^{-1}\psi(\eta,\xi),\qquad \psi(\eta,\xi) = \frac{1}{1+\eta\kappa(\xi)}.
\]
}

Following the ideas in \cite{2018AdjeridBenromdhaneLin}, \xz{the extended jump conditions \eqref{eq: jumpEx}, rewritten in Frenet coordinates as \eqref{eq: Fr jump3}}, lead to the following weak jump conditions on $\hGamma_{K_F}$:
\begin{eqnarray}
\int_{\hGamma_{K_F}} \hv \jump{ \hbeta \partial_\eta^j \mathscr{L}(\hu) }_{{\hat \Gamma}_{K_F}} d\xi = 0, \quad\forall \hv \in \mathbb{P}_{m-2 - j}(\hGamma_{K_F}), 
\label{eq: Fr jump3w}
\end{eqnarray}
for $j = 0, 1, \cdots, m-2.$

On each Frenet reference element $\hat{K}_F$ we construct the piecewise polynomial space 
\begin{equation}\label{eq: ref space}
\hat{\mathcal{V}}^m_{\hat\beta}(\hat{K}_F) = \left\{\hat \phi: \hat{K}_F\to \mathbb{R}: ~\hat{\phi}|_{\hat{K}_F^\pm}\in \mathbb{P}_m, ~\hat{\phi} \text{ satisfies \eqref{eq: Fr jump1},\eqref{eq: Fr jump2}, and \eqref{eq: Fr jump3w}}\right\}.
\end{equation}
where $\mathbb{P}_m$ is the polynomial space of order no more than $m$. The local geometric-conforming IFE space on the original interface triangle $K$ is defined by
\begin{equation}\label{eq: local space}
{\mathcal{V}}^m_{\beta}({K}) = \left\{\hat \phi\circ R_\Gamma|_K: ~\hat{\phi}\in \hat{\mathcal{V}}^m_{\hat\beta}(\hat{K}_F) \right\}.
\end{equation}
Then the global $\mathbb{P}_m$ GC-IFE space is defined as
\begin{equation}\label{eq: global space}
   {\mathcal{V}}^m_{\beta}({\mathcal{T}_h}) = \left\{v\in L^2(\Omega): v|_K\in {\mathcal{V}}^m_{\beta}({K}), \, \text{if}\, K\in{\mathcal{T}_h^i}; v\in \mathbb{P}_m(K) \, \text{if}\,  K\in{\mathcal{T}_h^n}
   \right\}. 
\end{equation}

In the following three sections, we \xz{introduce} several strategies to construct $\mathbb{P}_m$ GC-IFE space ${\mathcal{V}}^m_{\beta}({K})$ on interface triangles.

\section{Initial Construction of IFE Spaces on Triangles}
\label{sec: initial construction}
The first strategy to build $\mathbb{P}_m$ Frenet-IFE spaces on triangles follows the initial construction of the $\mathbb{Q}_m$ Frenet-IFE spaces on rectangular elements in \cite{2024AdjeridLinMeghaichi}.  We begin by constructing a monomial basis in $\mathbb{P}_m$ in Frenet coordinates $\eta$-$\xi$ on $\hK_F$. Then, through a change of variable with the Frenet transform, we obtain a set of IFE basis functions on fictitious element $K_F$. Finally, the Frenet IFE functions based on $\mathbb{P}_m$ 
polynomials are obtained by restricting these functions on $K_F$ to the triangular interface element $K$. 


Let $\left\{p_i(\xi)\right\}_{i = 0}^m$ be a basis of $\mathbb{P}_m(\widehat{\Gamma}_{K_F})$ with 
$\text{deg}(p_i) = i, ~0 \leq i \leq m$. It is straightforward to verify that the following polynomials form a basis for $\mathbb{P}_m(\hK_F)$:

\begin{eqnarray}
\begin{array}{ccccccccc}
p_0(\xi) &                             \\
p_1(\xi) & \eta p_0(\xi)      &        \\
p_2(\xi) & \eta p_1(\xi)      & \eta^2p_0(\xi) \\
p_3(\xi) & \eta p_2(\xi)      & \eta^2p_1(\xi) & \eta^3p_0(\xi) \\
\vdots   & \vdots             & \vdots         &\vdots  & \ddots \\
p_m(\xi) & \eta p_{m-1}(\xi)  & \eta^2 p_{m-2}(\xi)&\eta^3 p_{m-3}(\xi)  & \cdots & \eta^mp_0(\xi).
\end{array} \label{eq:basis_of_Pm}
\end{eqnarray}
There are $l_m = {(m-1)m}/{2}$ polynomials in columns $3$ through $(m+1)$ of the table \eqref{eq:basis_of_Pm}.
For convenience, we denote these basis polynomials by $\cN_i(\eta,\xi)$ with $i=1,2,\cdots,l_m$ where
\begin{eqnarray*}
\begin{array}{ccccc}
\cN_1  = \eta^2p_0(\xi) \\
\cN_2  = \eta^2p_1(\xi)          & \cN_{3} = \eta^3p_0(\xi) \\
\vdots          &\vdots & \ddots \\
\cN_{l_m-m+2}  = \eta^2 p_{m-2}(\xi) & \cN_{l_m-m+3} = \eta^3 p_{m-3}(\xi) & \cdots & \cN_{l_m} = \eta^mp_0(\xi)
\end{array}
\end{eqnarray*}
Following the same argument used in proving Lemma 3 of \cite{2024AdjeridLinMeghaichi} we obtain the analogous result for $\mathbb{P}_m$ polynomial spaces.
\begin{lemma}\label{lem:simple_ShFun_FrenetIFE_tri}
The piecewise polynomials
\begin{eqnarray}\label{eq: first several basis}
\hphi_{i, j}(\eta, \xi) = \frac{1}{\hbeta} \eta^i p_j(\xi), \quad 1 \leq i \leq m, \quad 0 \leq j \leq m - i
\label{eq:simple_ShFun_FrenetIFE_tri}
\end{eqnarray}
satisfy the interface jump conditions \eqref{eq: Fr jump1}, \eqref{eq: Fr jump2}, and  \eqref{eq: Fr jump3w}.
\end{lemma}

As stated in Lemma \ref{lem:simple_ShFun_FrenetIFE_tri}, the functions $\hphi_{i, j}$ in \eqref{eq:simple_ShFun_FrenetIFE_tri} provide $m(m+1)/2$ IFE shape functions on the Frenet fictitious element
$\hK_F$. We now describe the construction of the remaining $m+1$ IFE shape functions.

Using the basis in \eqref{eq:basis_of_Pm}, we can expand $\hphi^s(\eta, \xi) \in \mathbb{P}_m(\hK_F^s), ~s = \pm$ as follows:
\begin{equation*}
\hphi^s(\eta, \xi) = \sum_{k = 0}^m \left(\sum_{l = 0}^{m-k} C_{k, l}^s \eta^k p_l(\xi) \right) =
\sum_{k = 0}^m \eta^k \left(\sum_{l = 0}^{m-k} C_{k, l}^s p_l(\xi) \right) = \sum_{k = 0}^m \eta^k p_k^s(\xi), 
\end{equation*}
where 
\[
p_k^s(\xi) = \sum_{l = 0}^{m-k} C_{k, l}^s p_l(\xi) \in \mathbb{P}_{m-k}(\hGamma_{K_F}).
\]
Following the idea of Lemma 1 in \cite{2024AdjeridLinMeghaichi}, we obtain the unisolvency property.
\begin{lemma} \label{lem:FrenetIFE_extension_tri}
Given a polynomial $\hphi^s(\eta, \xi) \in \mathbb{P}_m(\hK_F^s)$, there exists a unique
$\hphi^{s'}(\eta, \xi) \in \mathbb{P}_m(\hK_F^{s'})$ such that \xz{the piecewise polynomial $\hat\phi$ defined by}
\begin{eqnarray*}
\hphi(\eta, \xi) = \begin{cases}
\hphi^-(\eta, \xi), & (\eta, \xi) \in \hK_F^-, \\
\hphi^+(\eta, \xi), & (\eta, \xi) \in \hK_F^+.
\end{cases}
\end{eqnarray*}
satisfies the jump conditions \eqref{eq: Fr jump1}, \eqref{eq: Fr jump2} and
\eqref{eq: Fr jump3w}.
\end{lemma}

For $i = 0, 1, \cdots, m$, Lemma \ref{lem:FrenetIFE_extension_tri} implies that an IFE shape function can be constructed in the form
\begin{equation}\label{eq:hphi_0i}
\hphi_{0, i}(\eta, \xi) = \begin{cases}
\hphi_{0, i}^-(\eta, \xi) = p_i(\xi), & (\eta, \xi) \in \hK_K^-, \\[5pt]
\hphi_{0, i}^+(\eta, \xi) = \sum\limits_{k = 0}^m \eta^k \left(\sum\limits_{l = 0}^{m-k} C_{k, l}^+ p_l(\xi) \right), & (\eta, \xi) \in \hK_K^+.
\end{cases}
\end{equation}
where coefficients $C_{k, l}^+$, $0 \leq l \leq m-k, ~0 \leq k \leq m$ are chosen so that
$\hphi_{0,i}(\eta, \xi)$ satisfies the interface jump conditions \eqref{eq: Fr jump1}, \eqref{eq: Fr jump2}, and  \eqref{eq: Fr jump3w}.
The jump condition \eqref{eq: Fr jump1} implies that
\begin{eqnarray*}
p_i(\xi) = \hphi_{0,i}^-(0, \xi) = \hphi_{0,i}^+(0, \xi) = \sum\limits_{l = 0}^{m} C_{0, l}^+ p_l(\xi), ~\forall \xi \in \widehat{\Gamma}_{K_F}.
\end{eqnarray*}
Hence, we should have
\begin{eqnarray*}
C_{0, l}^+ = \begin{cases}
1, & l = i, \\
0, & l \not= i
\end{cases}~~0 \leq l \leq m
\end{eqnarray*}
and
\begin{align*}
\hphi_{0,i}^+(\eta, \xi) &= p_i(\xi) + \sum_{k=1}^m \eta^k \left(\sum_{l=0}^{m-k} C_{k, l}^+ p_l(\xi)\right), \\
\frac{\partial \hphi_{0,i}^+(\eta, \xi)}{\partial \eta} &= \sum_{k=1}^m k\eta^{k-1} \left(\sum_{l=0}^{m-k} C_{k, l}^+ p_l(\xi)\right).
\end{align*}
%
Then, the jump condition \eqref{eq: Fr jump2} implies that
\begin{eqnarray*}
0 = \hbeta^- \frac{\partial \hphi_{0,i}^-(0, \xi)}{\partial \eta} = \hbeta^+ \frac{\partial \hphi_{0,i}^+(0, \xi)}{\partial \eta} = \sum\limits_{l = 0}^{m-1} C_{1, l}^+ p_l(\xi).
\end{eqnarray*}
Hence,
$C_{1, l}^+  = 0, ~0 \leq l \leq m-1,$
and
\begin{equation}\label{eq:hphi_0i+}
\begin{split}
\hphi_{0,i}^+(\eta, \xi) &= p_i(\xi) + \sum_{k=2}^m \eta^k \left(\sum_{l=0}^{m-k} C_{k, l}^+ p_l(\xi)\right) \\
&= p_i(\xi) + \sum_{l=1}^{l_m} C_l^i \cN_l(\eta, \xi), \quad 0 \leq i \leq m.
\end{split}
\end{equation}
We now use the jump conditions \eqref{eq: Fr jump3w} to derive the formulas for computing the coefficient $C_l^i, ~l = 1, 2, \cdots, \xz{l_m}$. 
These can be uniformly written as follows:
\begin{equation}\label{eq:coef_of_hphi_0i}
\begin{split}
\sum_{l=1}^{l_m} \left(\int_{\hGamma_{K_F}} p_k(\xi) \left(\frac{\partial^{j}}{\partial \eta^{j}}\mathscr{L}\left(\cN_l(0, \xi)\right)\right) d\xi\right) C_l^i \\
= \frac{\hbeta^- - \hbeta^+}{\hbeta^+} \int_{\hGamma_{K_F}} p_k(\xi) \left(\frac{\partial^{j}}{\partial \eta^{j}}\mathscr{L}(p_i(\xi))\right) d\xi
\end{split}
\end{equation}
for $~~0 \leq k \leq m-2-j,~0 \leq j \leq m-2.$

The linear system of ${l_m} = (m-1)m/2$  equations \eqref{eq:coef_of_hphi_0i} about the coefficients $C_l^i, ~~l = 1, 2, \cdots, {l_m} $ can be written as 
\begin{equation}\label{eq: system}
    A\mathbf{c}^{(i)} =\frac{\hbeta^- - \hbeta^+}{\hbeta^+}\mathbf{b}{(i)}
\end{equation}
where the matrix and right-hand side are defined by
\[
A = \begin{bmatrix}
    A^{(0)}\\A^{(1)}\\\vdots\\A^{(m-2)}
\end{bmatrix},\qquad
\mathbf{b}(i) = \begin{bmatrix}
    \mathbf{b}^{(0)}(i)\\\mathbf{b}^{(1)}(i)\\\vdots\\\mathbf{b}^{(m-2)}(i)
\end{bmatrix},\qquad
\]
and
\begin{align*}
A_{k,l}^{(j)} &= \int_{\hGamma_{K_F}} p_k(\xi) \left(\frac{\partial^{j}}{\partial \eta^{j}}\mathscr{L}\left(\cN_l(0, \xi)\right)\right) d\xi, \\
\mathbf{b}_k^{(j)}(i) &= \int_{\hGamma_{K_F}} p_k(\xi) \left(\frac{\partial^{j}}{\partial \eta^{j}}\mathscr{L}(p_i(\xi))\right) d\xi.
\end{align*}
The system \eqref{eq: system} is uniquely solvable, as guaranteed by Lemma \ref{lem:FrenetIFE_extension_tri}. Once the remaining $m+1$ basis functions $\hphi_{0, i}$ are determined, together with $\hphi_{i,j}$ in \eqref{eq: first several basis}, we obtain all $(m+1)(m+2)/2$ basis functions. The local space  is therefore
\begin{equation}\label{eq: ref space 1}
\hat{\mathcal{V}}^m_{\hat\beta}(\hat{K}_F) = \text{span} \left\{\hat \phi_{i,j}: 0\le i\le m, 0\le j\le i \right\}.
\end{equation}

\begin{remark}
The construction developed above closely follows the construction of Frenet IFE spaces on rectangles presented in \cite{2024AdjeridLinMeghaichi}. A key advantage of this approach is that the majority of basis functions can be constructed explicitly via \eqref{eq:simple_ShFun_FrenetIFE_tri}. For a polynomial space of degree $m$, $m(m+1)/2$ of the $(m+1)(m+2)/2$ basis functions can be written  explicitly, and only the remaining $m+1$ basis functions are obtained by solving a $l_m\times l_m$ linear system \eqref{eq:coef_of_hphi_0i}. 
\end{remark}

\begin{remark}
Moreover, letting 
\begin{eqnarray*}
{\hat \phi}_1(\eta, \xi) = p_0(\xi), ~{\hat \phi}_2(\eta, \xi)= p_1(\xi), ~{\hat \phi}_3(\eta, \xi)= \eta p_0(\xi)/{\hat \beta}(\eta, \xi),
\end{eqnarray*}
we can easily show that  
\begin{equation}\label{eq: ref space 1_m=1}
\hat{\mathcal{V}}^1_{\hat\beta}(\hat{K}_F) = \text{span} \left\{{\hat \phi}_1, ~{\hat \phi}_2, ~{\hat \phi}_3\right\}.
\end{equation}
Consequently, the local $\mathbb{P}_1$ Frenet IFE space on an interface element $K$ has the following explicit expression, no numerical calculation is required.  
\begin{equation}\label{eq: local space_m=1}
{\mathcal{V}}^1_{\beta}({K}) = \left\{{\hat \phi}_1\circ R_\Gamma|_K, ~{\hat \phi}_2\circ R_\Gamma|_K, ~{\hat \phi}_3\circ R_\Gamma|_K \right\}.
\end{equation}
This is in contrast to the $\mathbb{P}_1$ IFE space \cite{2004LiLinLinRogers} whose basis functions have to numerically determined on every interface element. 

\end{remark}

\begin{figure}[h]
\includegraphics[width = .47\textwidth]{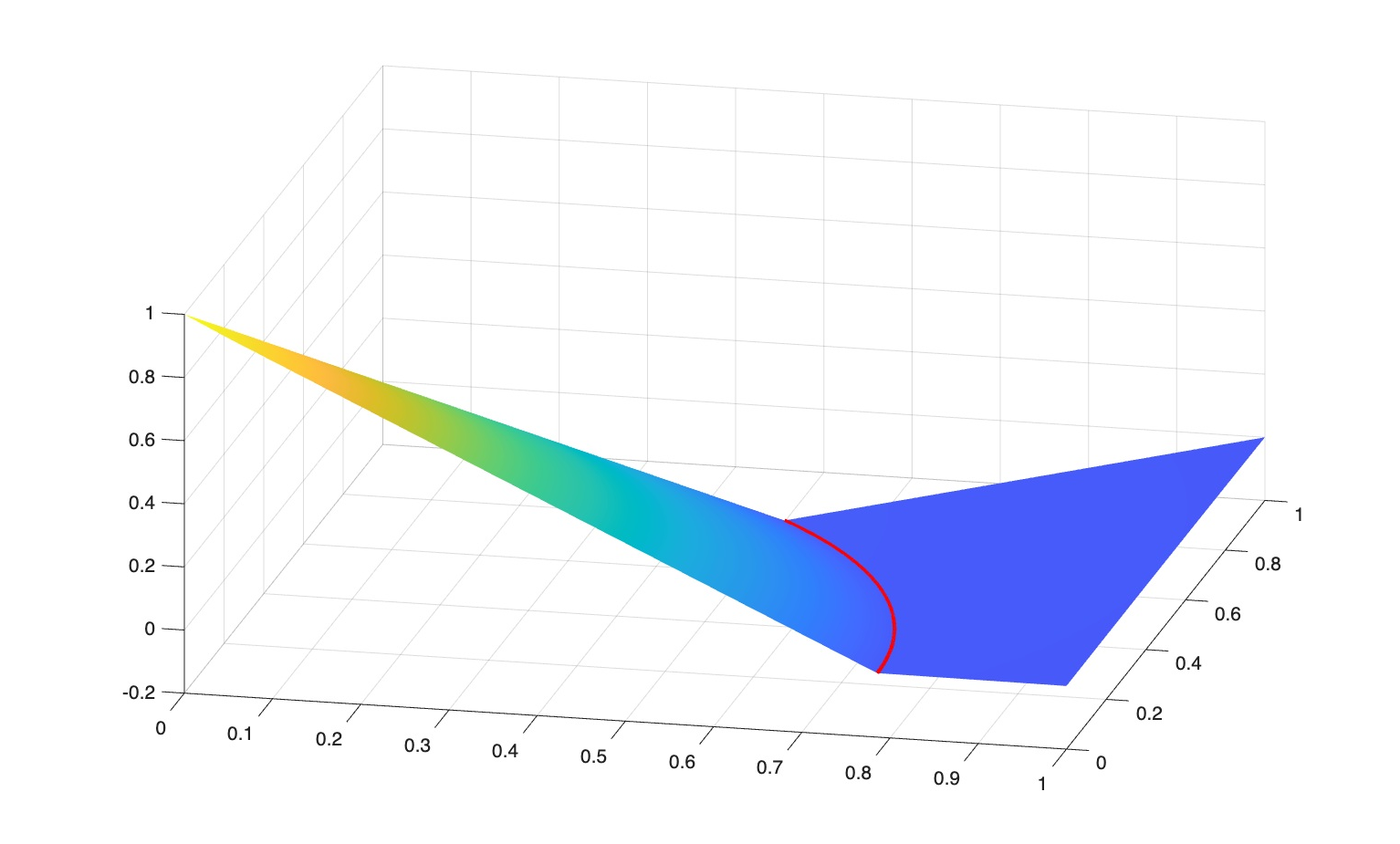}
\includegraphics[width = .47\textwidth]{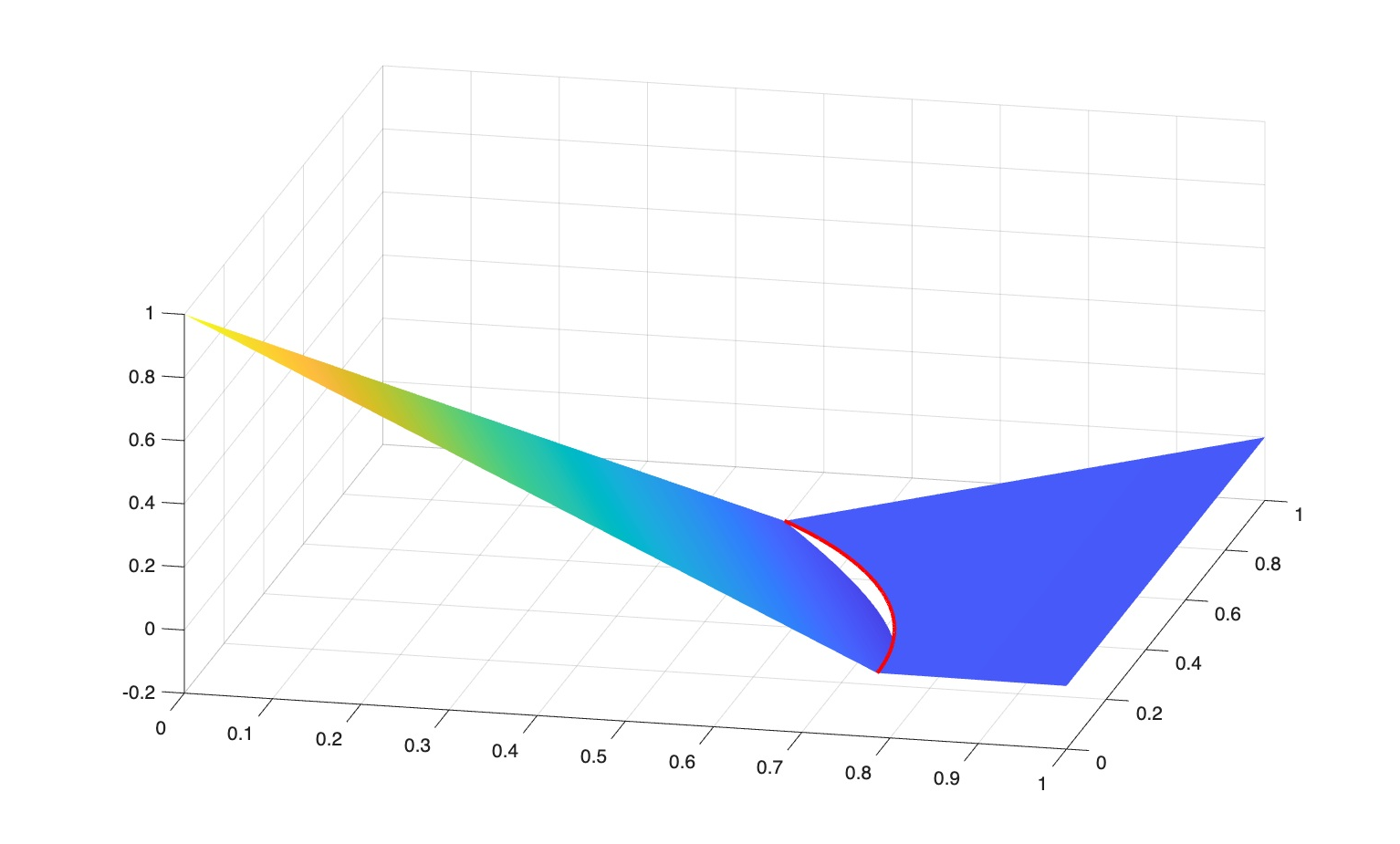}
\caption{A comparison of Frenet-IFE shape function (left) and a classical $P_1$ IFE shape function (right)}.
\label{fig: FIFE-P1}
\end{figure}

\begin{remark}
The Frenet–IFE shape functions $\phi$ are geometrically conforming: the pieces $\phi^+$ and $\phi^-$ match along the entire interface curve within the element, as shown in the left plot of Figure~\ref{fig: FIFE-P1}. In contrast, classical $\mathbb{P}_1$ IFE shape functions \cite{2004LiLinLinRogers} enforce continuity only at the two intersection points where the interface cuts the element edges. 

\end{remark}

\section{General Construction for Frenet-IFE Basis on Triangles}
\label{sec: general construction basis}
It is known that monomial bases tend to produce ill-conditioned linear systems, particularly when the polynomial degree becomes large. For example, using the basis $\{x^i\}_{i=0}^m$ to compute the $L^2$ projection onto $\mathbb{P}_m$ leads to a Hilbert matrix, whose condition number grows rapidly and is notoriously unfavorable for numerical computation. In contrast, orthogonal polynomial bases, such as Legendre polynomials, often yield much better conditioned systems due to their intrinsic orthogonality and reduced numerical coupling among the basis functions. These considerations motivate the use of more general and computationally robust bases when constructing the GC-IFE space $\mathcal{V}^m_\beta(K)$.

\subsection{General Construction of GC-IFE Basis}

\xz{
Following the construction in \cite{2025AdjeridLinMeghaichi2}, we define a set of polynomials $\{R_l: 1\le l\le (m+1)(m+2)/2\}$ as follows
\begin{equation}
R_{(t+s)(t+s+1)/2 + t+1}(\eta, \xi) =  q_t
\left(\frac{\eta}{\eta_h}\right)p_s\left(\frac{\xi-\bar\xi}{\xi_h}\right), ~0 \leq t \leq m, 0 \leq s\leq m-t 
\end{equation}
where $q_0 = p_0\equiv 1$, and $\text{deg}(q_t) = t$, $\text{deg}(p_s) = s$ when $s,t\ge 1$ and such that
\[
q_1(0) = 0, \quad q_i(0) = q_i'(0) =0, ~\quad 2 \leq i \leq m.
\]
It is straightforward to verify that Let $\{R_l\}_{l=1}^{d_m}$ form a polynomial basis for the space $\mathbb{P}_m(K_F)$ in the ($\eta$,$\xi$)-coordinates where $d_m = (m+1)(m+2)/2$. 
}
We construct a set of basis functions $\{\hat\lambda_j: 1\le j\le d_m\}$  on the Frenet reference element $\hat K_F$ of the following form:
\begin{equation}\label{eq: basis general}
\hat\lambda_j(\eta,\xi) = 
\begin{cases}
\hat\lambda_j^-(\eta,\xi) = \displaystyle\sum\limits_{i=1}^{d_m} C_{ij}^-R_i(\eta,\xi), & \eta <0, \\
\hat\lambda_j^+(\eta,\xi) = \displaystyle\sum\limits_{i=1}^{d_m} C_{ij}^+R_i(\eta,\xi), &\eta >0,
\end{cases} \qquad\  \xz{1} \leq j \leq d_m.
\end{equation}
The coefficients $C_{ij}^\pm$  are to be determined so that $\hat\lambda_j$ satisfies the interface jump conditions on ${\widehat \Gamma}_{K_F}$. These conditions give rise to the following system of equations:
\begin{equation}\label{eq: sys1}
\begin{split}
&\sum_{i=1}^{d_m} \left(\int_{\xi_{0,K}}^{\xi_{1,K}} R_i(0,\xi) \xz{p_k(\xi)} \, d\xi\right) C_{ij}^+ \\
=& \sum_{i=1}^{d_m} \left(\int_{\xi_{0,K}}^{\xi_{1,K}} R_i(0,\xi)\xz{p_k(\xi)} \, d\xi\right) C_{ij}^-, \quad k=0,1,\ldots,m.
\end{split}
\end{equation}

\begin{equation}\label{eq: sys2}
\begin{split}
&\sum_{i=1}^{d_m} \left(\int_{\xi_{0,K}}^{\xi_{1,K}} \partial_\eta R_i(0,\xi)\xz{p_k(\xi)} \, d\xi\right) C_{ij}^+ \\
=& \frac{\beta^-}{\beta^+} \sum_{i=1}^{d_m} \left(\int_{\xi_{0,K}}^{\xi_{1,K}} \partial_\eta R_i(0,\xi)\xz{p_k(\xi)} \, d\xi\right) C_{ij}^-, \quad k=0,\ldots,m-1.
\end{split}
\end{equation}

\begin{equation}\label{eq: sys3}
\begin{split}
&\sum_{i=1}^{d_m} \left(\int_{\xi_{0,K}}^{\xi_{1,K}} \partial_\eta^n \mathscr{L}(R_i(0,\xi))\xz{p_k(\xi)} \, d\xi\right) C_{ij}^+ \\
=& \frac{\beta^-}{\beta^+} \sum_{i=1}^{d_m} \left(\int_{\xi_{0,K}}^{\xi_{1,K}} \partial_\eta^n \mathscr{L}(R_i(0,\xi))\xz{p_k(\xi)} \, d\xi\right) C_{ij}^-, \\
&\qquad\qquad k=0,\ldots,m-2-n, \quad n=0,\ldots,m-2.
\end{split}
\end{equation}

Let $C_j^\pm$ denote the column vector consisting of $C^\pm_{ij}$ for $1\le i\le d_m$. Then the linear system \eqref{eq: sys1}-\eqref{eq: sys3} can be written in the following matrix form 
\begin{equation}\label{eq: coef}
\tilde A C_j^+ = J \tilde A C_j^-
\end{equation}
where \xz{$\tilde A$ is the coefficient matrix associated with \eqref{eq: sys1}-\eqref{eq: sys3}}, and  
\[J = \text{diag}(1,1,\cdots,1, \frac{\beta^-}{\beta^+},\cdots,\frac{\beta^-}{\beta^+})\] 
is a diagonal matrix whose first $m+1$ entries equal $1$ and the remaining $(d_m - m - 1)$ entries equal $\beta^- / \beta^+$. 

The existence and uniqueness of GC-IFE basis functions imply that, for each fixed $j$, the vector $C_j^+$ is uniquely determined from \eqref{eq: coef} once $C_j^-$ is given. 
Thus, given any basis ${\hat\lambda_j^- : 1\le j\le d_m}$ of $\mathbb{P}_m$ on $\hat K_F^-$, the extension mapping \eqref{eq: coef} uniquely determines the corresponding functions
${\hat\lambda_j^+ : 1\le j\le d_m}$ on $\hat K_F^+$. 
In this way, the complete set of local GC-IFE basis functions ${\hat\lambda_j(\eta,\xi) }$, $ 1 \le j \le d_m$ is generated.

Let $C^\pm = [C_1^\pm, C_2^\pm, \dots, C_{d_m}^\pm]$ denote the coefficient matrices. Then \eqref{eq: coef} is equivalent to the matrix identity
\begin{equation}\label{eq: coef matrix}
\tilde A C^+ = J \tilde A C^-.
\end{equation}
In practice, one may select any nonsingular matrix $C^-$ and compute $C^+$ directly from \eqref{eq: coef matrix}; or conversely, prescribe $C^+$ and compute $C^-$.

Using the general construction, the local space can be written as 
\begin{equation}\label{eq: ref space 1}
\hat{\mathcal{V}}^m_{\hat\beta}(\hat{K}_F) = \text{span} \left\{\hat \lambda_{j}:  1\le j\le d_m \right\}.
\end{equation}

\begin{remark}
As discussed in \cite{2025AdjeridLinMeghaichi2}, the initial construction of the Frenet-based IFE basis functions $\phi_{ij}$ presented in Section~\ref{sec: initial construction} arises as a special case of this more general GC-IFE framework. The formulation introduced here therefore not only unifies the two constructions but also provides greater flexibility in selecting polynomial bases that may improve conditioning and numerical performance.

\end{remark}

\subsection{Efficient generation of matrices and vectors}
Next we present an approach to efficiently generate local mass and stiffness matrices and source vectors on a triangular interface element using GC-IFE basis functions. Let $\{\hat\lambda_j(\eta,\xi): j=1,2,\cdots, d_m\}$ be a basis for the reference Frenet space $\hat{\mathcal{V}}^m_{\hat\beta}(\hat{K}_F)$ via general construction. Then $\{\lambda_j:=\hat\lambda_j\circ R_\Gamma, j=1,2,\cdots, d_m\}$ is a basis for the local GC-IFE space $V^m(K)$ on an interface element $K$. The local mass matrix $M_K = (m_{ij})$ can be written as 
\begin{equation}\label{eq: integral}
m_{ij} = \int_{K}  \lambda_i\lambda_j d\mathbf{x} = \int_{K}  (\hat\lambda_i\circ R_\Gamma)(\hat\lambda_j\circ R_\Gamma) d\mathbf{x} 
=\sum_{s=\pm} \int_{K^s}  (\hat\lambda_i\circ R_\Gamma)(\hat\lambda_j\circ R_\Gamma) d\mathbf{x}.
\end{equation}
On each piece ${K^s}$, the integration \eqref{eq: integral} is approximated by numerical quadratures, i.e.,
\begin{equation}
m_{ij} \approx \sum_{s=\pm} \sum_{k=1}^{n_q^s} w_k^s \big(\hat\lambda_i\circ R_\Gamma(\mathbf{x}_k^s)\big)\big(\hat\lambda_j\circ R_\Gamma({\mathbf{x}}_k^s)\big) = \sum_{s=\pm} \sum_{k=1}^{n_q^s} w_k^s \hat\lambda_i(\hat{\mathbf{x}}_k^s)\hat\lambda_j(\hat{\mathbf{x}}_k^s),
\end{equation}
where $w_k^s$ and $\mathbf{x}_k^s$ are quadrature weights and nodes on $K^s$ and $n_q^s$ is the number of quadrature points on $K^s$ in a certain quadrature rule. Note that $K^s$ is either a curved triangle or quadrilateral. The Stroud quadrature rule introduced in Section 3.2 of \cite{2025AdjeridLinMeghaichi2}  for rectangular meshes \xz{is also} applicable in this case. 

Denote the Vandermonde matrices associated with the basis $\hat\lambda_i$ and the basis $R_i$ of $\mathbb{P}_m$ to be 
\begin{equation}
V^s = \left(\hat\lambda_i(\hat{\mathbf{x}}_k^s)\right)_{k=1,i=1}^{n_q^s,d_m}
,\qquad
L^s = \left(R_i(\hat{\mathbf{x}}_k^s)\right)_{k=1,i=1}^{n_q^s,d_m},\quad s = \pm.
\end{equation}
Then the relation \eqref{eq: basis general} implies $V^s = L^sC^s$. The mass matrix $M_K$ on the interface element $K$ can be approximated by $M_{K,q}$ of the following form
\begin{equation}\label{eq:mass}
\begin{split}
\xz{M_K \approx M_{K,q}} &= \left(\sum_{s=\pm} \sum_{k=1}^{n_q^s} w_k^s \hat\lambda_i(\hat{\mathbf{x}}_k^s)\hat\lambda_j(\hat{\mathbf{x}}_k^s)\right)_{i,j=1}^{d_m} \\
&= \sum_{s=\pm} (V^s)^T W^s V^s = \sum_{s=\pm} (C^s)^T(L^s)^T W^s L^sC^s.
\end{split}
\end{equation}
where $W^s = \text{diag}(w_1^s,w_2^s,\cdots,w_{n_q^s}^s)$ is the diagonal matrix consisting of quadrature weights. 

\begin{remark}
The notation $M_{K,q}$ emphasizes that the approximation depends on the
quadrature rule used for integration.
\end{remark}
\begin{remark}
The procedure in \eqref{eq:mass} yields an efficient approximation of the element mass matrix $M_K$, because the matrices $L^s$ are assembled by evaluating the standard polynomial basis
functions $R_i$ at the quadrature points $\hat{\mathbf{x}}_k^s$, rather than using piecewise polynomials.  
After the coefficients of the GC-IFE shape functions $C^\pm$ are computed,
the mass matrix $M_K$ is obtained by standard matrix multiplication.
\end{remark}

The construction of the local right hand side vector $\mathbf{f}_K$ on an interface element follows a similar procedure. In fact, 
\begin{equation}
{\bf f}_K \approx {\bf f}_{K,q}  = \sum_{s=\pm} (V^s)^TW^s\mathbf{r}^s =  \sum_{s=\pm} (C^s)^T(L^s)^TW^s\mathbf{r}^s,\qquad \end{equation}
\text{where}
\[
\mathbf{r}^s = 
\Big(f(\mathbf{x}_1^s),f(\mathbf{x}_2^s),\cdots,f(\mathbf{x}_{n_q^s}^s)\Big)^T.
\]

This approach is applicable to the stiffness matrix $S$ as well. 
The entries of $\xz{S_K} = (S_{ij})$ can be written as 
\begin{equation}\label{eq: stiff local}
S_{ij} = 
\int_{K}  \beta \nabla  \lambda_{i}\cdot \nabla \lambda_{j} d\mathbf{x} 
=\sum_{s=\pm} \int_{K^s}  \beta^s \nabla (\hat\lambda_i\circ R_\Gamma)\cdot  \nabla (\hat\lambda_j\circ R_\Gamma) d\mathbf{x}.
\end{equation}
Note that 
\begin{equation}\label{eq:stiffness}
\begin{split}
S_{ij} &\approx \sum_{s=\pm} \sum_{k=1}^{n_q^s} w_k^s \beta^s \nabla(\hat\lambda_i \circ R_\Gamma(\mathbf{x}_k^s)) \cdot \nabla(\hat\lambda_j \circ R_\Gamma(\mathbf{x}_k^s)) \\
&= \sum_{s=\pm} \sum_{k=1}^{n_q^s} w_k^s \beta^s \big(\hat\lambda_{i,\eta}(\hat x_k^s)\hat\lambda_{j,\eta}(\hat x_k^s) + \rho^2 \hat\lambda_{i,\xi}(\hat x_k^s) \hat\lambda_{j,\xi}(\hat x_k^s)\big).
\end{split}
\end{equation}
Denote the Vandermonde matrices:
\[
V_\eta^s = \left(\hat\lambda_{i,\eta}(\hat{\mathbf{x}}_k^s)\right)_{k=1,i=1}^{n_q^s,d_m},\qquad
V_\xi^s = \left(\hat\lambda_{i,\xi}(\hat{\mathbf{x}}_k^s)\right)_{k=1,i=1}^{n_q^s,d_m}
\]
and 
\[
L_\eta^s = \left(R_{i,\eta}(\hat{\mathbf{x}}_k^s)\right)_{k=1,i=1}^{n_q^s,d_m},\qquad 
L_\xi^s = \left(R_{i,\xi}(\hat{\mathbf{x}}_k^s)\right)_{k=1,i=1}^{n_q^s,d_m},\quad s = \pm.
\]

Then the stiffness matrix $S_K$ can be approximated by $S_{K,q}$ which can be written as 
\begin{equation}\label{eq:Sq}
\begin{split}
\xz{S_{K,q}} &= \sum_{s=\pm} \Big(\beta^s (V_\eta^s)^T W^s V_\eta^s + \beta^s \rho^2 (V_\xi^s)^T W^s V_\xi^s\Big) \\
&= \sum_{s=\pm} \Big(\beta^s (C^s)^T (L_\eta^s)^T W^s L_\eta^s C^s + \beta^s \rho^2 (C^s)^T (L_\xi^s)^T W^s L_\xi^s C^s\Big).
\end{split}
\end{equation}

Once again, the generation of the stiffness matrix requires the evaluation of the first derivative of the polynomial basis $R_i$ at the quadrature points.

\section{Reconstruction of Frenet-IFE Basis on Triangles}\label{sec: reconstruction basis}

The GC-IFE bases in the general construction \eqref{eq: basis general} focus primarily on enforcing the jump conditions, but not on computational performance when used to approximate interface problems. In this section, we introduce two reconstruction approaches that significantly improve the conditioning of the GC-IFE mass matrix and enhance computational performance.
%
%

We write Frenet IFE basis function $\hat \lambda_j$ via general construction as follows
\[\hat \lambda_j(\eta,\xi; C_j^-,C_j^+), \quad j=1,\cdots, d_m.\]
We include $C_j^-,C_j^+$ in the basis functions to emphasize that these bases can be completely determined by $C_j^-,C_j^+$ once the regular basis $R_j$ for $\mathbb{P}_m$ is chosen. 

We aim to introduce a reconstructed Frenet IFE basis $\tilde \lambda$ so that the resulting mass matrix \xz{is optimally conditioned}. To start, we denote the basis of $\hat{V}^m(\hat K_F)$ by
\[
\hat{B}(C) = \hat{B}(C^-,C^+) := \{\hat \lambda_j(\eta,\xi; C_j^-,C_j^+): j=1,\cdots, d_m\}.
\]
We first recall the Lemma 1 from  \cite{2025AdjeridLinMeghaichi}:
\begin{lemma}
Let $\hat{B}(C)$ be a basis of $\hat{\mathcal{V}}^m_{\hat\beta}(\hat{K}_F)$ and let $Q$ be a nonsingular matrix of size $d_m\times d_m$, then $\hat{B}(CQ)$ is also a basis of $\hat{\mathcal{V}}^m_{\hat\beta}(\hat{K}_F)$.
\end{lemma}

Let $\tilde C = CQ$. By \eqref{eq:mass}, the mass matrix $M(\hat{B}(\tilde C))$ associated with the basis $\hat{B}(\tilde C) = \hat{B}(\tilde C^-,\tilde C^+) = \hat{B}(C^-Q, C^+Q) $ can be approximated by 
\begin{equation}\label{eq:Mq_change_of_basis}
\begin{split}
\xz{M_{K,q}}(\hat{B}(\tilde C)) &= \sum_{s=\pm} (\tilde C^s)^T (L^s)^T W^s L^s \tilde C^s \\
&= \sum_{s=\pm} (C^s Q)^T (L^s)^T W^s L^s C^s Q = Q^T \xz{M_{K,q}}(\hat{B}(C)) Q.
\end{split}
\end{equation}
Hence, we can find a new basis $\hat{B}(\tilde C)$ such that the mass matrix $\xz{M_{K,q}}(\hat{B}(\tilde C))$ is better conditioned than the mass matrix $\xz{M_{K,q}}(\hat{B}(C))$. Inspired by \cite{2025AdjeridLinMeghaichi2}, we consider the following two approaches.

\subsection*{Reconstruction Approach 1} 
We choose $Q$ such that the mass matrix becomes the identity matrix. Let $V_1\Lambda V_1^T = \xz{M_{K,q}}(\hat{B}(C))$ be the singular value decomposition of the matrix $\xz{M_{K,q}}(\hat{B}(C))$ where $\Lambda$ denote the diagonal matrix of positive eigenvalues of $\xz{M_{K,q}}(\hat{B}(C))$. Let $Q_1 = V_1\Lambda^{-1/2}$ and $C^{(1)} = CQ_1$. Then we have 
\begin{equation}\label{eq:Mq_identity}
\begin{split}
\xz{M_{K,q}}(\hat{B}(C^{(1)})) &= \xz{M_{K,q}}(\hat{B}(CQ_1)) 
= Q_1^T \xz{M_{K,q}}(\hat{B}(C)) Q_1 \\
&= \Lambda^{-1/2} V_1^T \xz{M_{K,q}}(\hat{B}(C)) V_1 \Lambda^{-1/2} = I_{d_m\times d_m}.
\end{split}
\end{equation}

This approach transforms the GC-IFE basis $\hat{{B}}(C)$ into an orthonormal basis $\hat{{B}}(C^{(1)})$ with respect to the $L^2$ inner product and the chosen quadrature rule. We denote this reconstructed basis by
\begin{equation}
\hat{B}(C^{(1)}) = \hat{B}(C^-Q_1,C^+Q_1) := \{\hat \lambda^{(1)}_j: j=1,\cdots, d_m\}.
\end{equation}
The resulting approximate mass matrix $\xz{M_{K,q}}(\hat{{B}}(C^{(1)}))$ achieves the optimal condition number in theory. In actual computation with finite precision, the conditioning could be affected slightly as we can see in next section. 

\subsection*{Reconstruction Approach 2} The second approach uses a factorization of $\xz{M_{K,q}}(\hat{B}(C))$ as follows. 
\xz{Let $n_q = n_q^- + n_q^+$, and combine the quadrature points and weights:}
\[
\big\{(\mathbf{x}_r, w_r)\big\}_{r=1}^{n_q}   = \big\{(\mathbf{x}_k^-, w_k^-)\big\}_{k=1}^{n_q^-}
      \,\cup\,     \big \{(\mathbf{x}_k^+, w_k^+)\big\}_{k=1}^{n_q^+},
\]
then, by \eqref{eq:mass} we can express $\xz{M_{K,q}}(\hat{B}(C))$ as
\begin{equation}
\xz{M_{K,q}}(\hat{B}(C))
  = V(\hat{B}(C))^T W V(\hat{B}(C)),
\end{equation}
with 
$
W = \operatorname{diag}(w_1, w_2, \ldots, w_{n_q}),
$
and $V(\hat{B}(C))$ is the generalized Vandermonde matrix of the basis functions $\hat{B}(C)$ in with entries
$\mathbf{x}_r$, $1 \le r \le n_q$, i.e.,
\[
V_{r,j}(\hat{B}(C))
 = \hat{\lambda}_j(\mathbf{x}_r),
 \qquad
 1 \le r \le n_q,\; 1 \le j \le d_m.
\]
Let $\tilde{V}(\hat{B}(C)) = W^{1/2}\, V(\hat{B}(C))$, then we can factorize the matrix $\xz{M_{K,q}}(\hat{B}(C))$ as follows
\[
\xz{M_{K,q}}(\hat{B}(C)) = \tilde{V}(\hat{B}(C))^T \tilde{V}(\hat{B}(C)).
\]
Assume that $n_q \ge d_m$, we can compute the reduced SVD 
\[
\tilde{V} = U_2 \Sigma V_2^T,
\]
where $U_2$ is an $n_q \times d_m$ matrix and 
$\Sigma, V_2$ are $d_m \times d_m$ matrices.
Let $Q_2 = V_2 \Sigma^{-1}$ and $C^{(2)} = CQ_2$, then we obtain the second reconstructed basis 
\begin{equation}
\hat{B}(C^{(2)}) = \hat{B}(C^-Q_2,\,C^+Q_2) := \{\hat \lambda^{(2)}_j: j=1,\cdots, d_m\}.
\end{equation}
The associated mass matrix is
\[
\begin{split}
\xz{M_{K,q}}(\hat{{B}}(C^- Q_2,\, C^+ Q_2))
& = Q_2^T \xz{M_{K,q}}(\hat{{B}}(C)) Q_2  \\
& = \Sigma^{-1} V_2^T \tilde{V}^T \tilde{V} V_2 \Sigma^{-1} 
 = I_{d_m\times d_m}.
\end{split}
\]
Both reconstruction approaches make approximate mass matrix optimally conditioned. Moreover, when the eigenvalues of $\xz{M_{K,q}}(\hat{{B}}(C))$ are
distinct and computations are exact, the two approaches are mathematically equivalent. In this case,
they generate the same GC-IFE basis, i.e.,
\[
\hat{B}(C^{(1)}) = \hat{B}(C^{(2)})
\]
up to signs of the corresponding basis functions.



\begin{remark} \label{rem: condition}
In finite-precision computations, the two reconstruction approaches behave slightly differently. The condition number of $\xz{M_{K,q}}(\hat{B}(C))$ is the square of the
condition number of $\tilde{V}(\hat{B}(C))$ since $\Lambda = \Sigma^2$. As a result,
$\xz{M_{K,q}}(\hat{B}(C))$ is typically  more ill-conditioned than $\tilde{V}$, and the numerical
accuracy of the SVD in reconstruction approach (and therefore of $Q_1$) is limited by round-off errors in the small
singular values of $\xz{M_{K,q}}(\hat{B}(C))$. Consequently, the matrix $Q_1$ obtained numerically may
deviate substantially from its exact counterpart.
In contrast, the SVD of $\tilde{V}(\hat{B}(C))$ in reconstruction approach~2 is less sensitive to round-off
errors because $\tilde{V}$ has a much smaller condition number. Thus the transformation matrix $Q_2$ is generally more accurate. These results are consistent with the reconstruction for GC-IFE basis for rectangular elements \cite{2025AdjeridLinMeghaichi2}. 
\end{remark}
\section{Numerical Experiments}\label{sec: numerical}
In this section, we present a set of numerical experiments to demonstrate the performance of the geometrically conforming immersed finite element spaces on triangular meshes.

\begin{figure}[h]
\includegraphics[width = .35\textwidth]{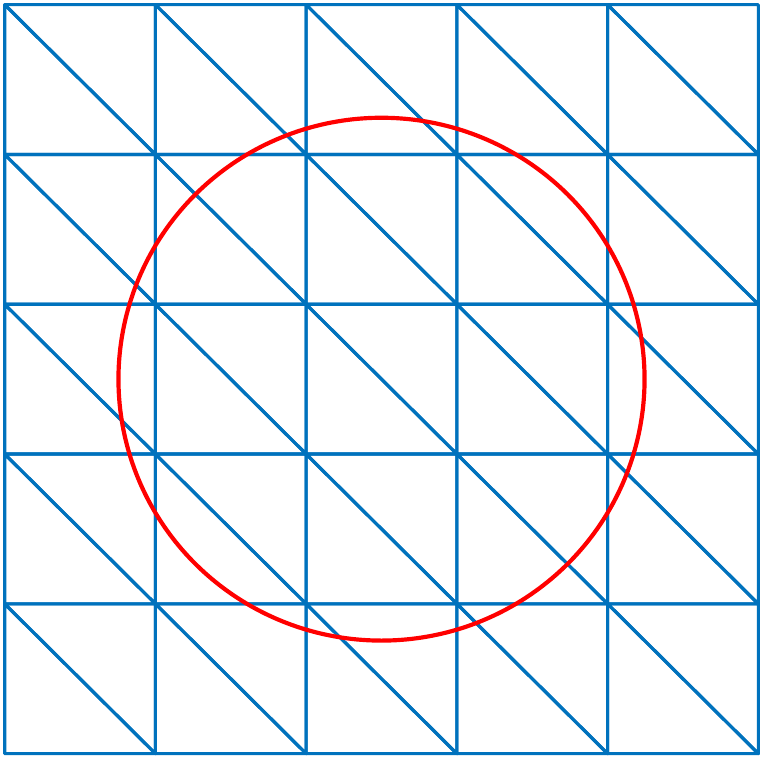}\qquad\quad
\includegraphics[width = .32\textwidth]{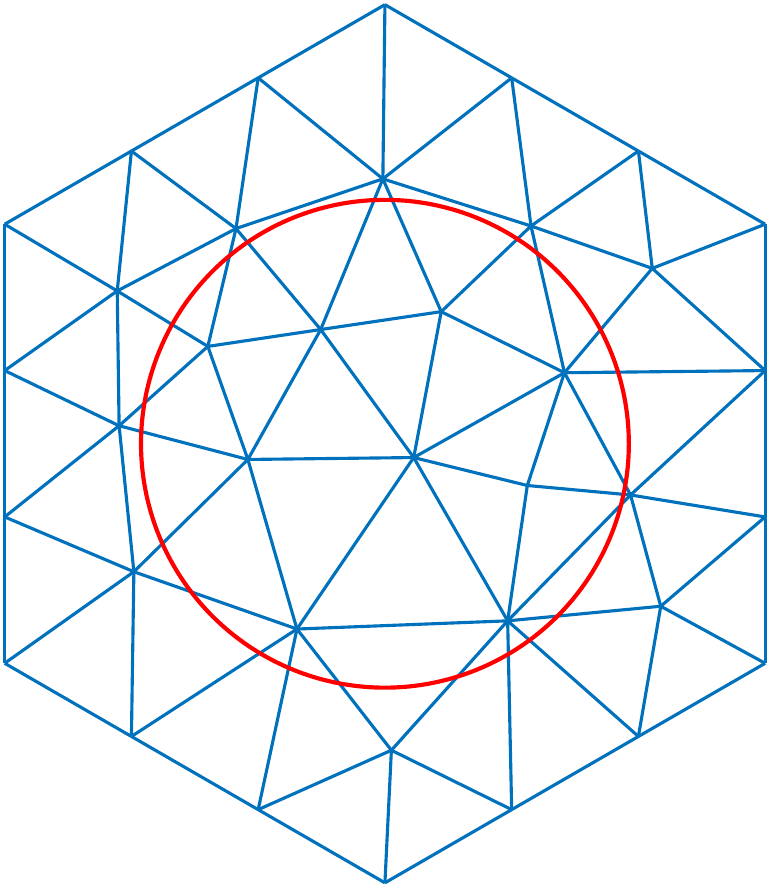}
\caption{A structured triangular mesh of a rectangular domain and an unstructured mesh of a hexagonal domain.}
\label{fig: mesh}
\end{figure}

\xz{We consider two types of computational domains. In the first three examples, the computational domain is a square and structured Cartesian triangular meshes are used. In the fourth example, we consider a hexagonal domain on which general unstructured meshes are used, see Figure \ref{fig: mesh}.}

In the first three examples, the computational domain is set to be a square $\Omega = (-1, 1)^2$ with a circular interface $\Gamma = \{ |\mathbf{x}| = r_0 \}$ where $r_0 = 1/\sqrt{3}$, which partitions $\Omega$ into two subdomains 
$\Omega^- = \{ |\mathbf{x}| < r_0 \}$ 
and $\Omega^+ = \{ |\mathbf{x}| > r_0 \}$.
The diffusion coefficient $\beta(x)$ is piecewise constant  $\beta(\mathbf{x})|_\Omega^\pm = \beta^\pm$ where $\beta^- = 1$ and $\beta^+$ is varied in the tests.

The exact solution is given by:
\begin{equation}\label{eq: true solution}
u(\mathbf{x}) =
\begin{cases}
	\dfrac{1}{\beta^-}\cos(\pi |\mathbf{x}|^2), & \mathbf{x} \in \Omega^-, \\[5pt]
	\dfrac{1}{\beta^+}\cos(\pi |\mathbf{x}|^2) + \cos(\pi r_0^2)\left(\dfrac{1}{\beta^-} - \dfrac{1}{\beta^+}\right), & \mathbf{x} \in \Omega^+,
\end{cases}
\end{equation}
which satisfies the interface jump conditions \eqref{eq: jump1}-\eqref{eq: jump2} as well as the Laplacian extended jump condition \eqref{eq: jumpEx}. We examine three types of basis:  $\{\hat\lambda_i\}_{i=1}^{d_m}$ from general construction, $\{\hat\lambda_i^{(1)}\}_{i=1}^{d_m}$ from the \xz{Reconstruction Approach 1}, and $\{\hat\lambda_i^{(2)}\}_{i=1}^{d_m}$ from \xz{Reconstruction Approach 2}. Polynomial degrees up to $m = 9$ are tested. 

Numerical tests have been carried out on both structured and unstructured triangular meshes; see Figure~\ref{fig: mesh} for representative examples. Since the results are comparable, we only report the data obtained from unstructured meshes, which are generated using the MATLAB PDE Toolbox.


\subsection{Conditioning of Mass and Stiffness Matrices}\label{sec: condition}

We first investigate the conditioning of the mass \xz{and stiffness} matrices associated with different choices of basis functions. Local mass matrices are computed on all interface elements of \xz{a structured $10\times 10$ triangular mesh}, and the maximum condition number over all interface elements is reported in Table~\ref{tab:condMlocal}. Three types of basis functions are considered. For the general construction, the condition numbers grow rapidly as both the polynomial degree and the jump ratio increase. In contrast, both reconstruction approaches appear \xz{less sensitive} to these variations. In particular, the Reconstruction Approach 2 demonstrates superior robustness at higher polynomial degrees. When the degree exceeds 7, the condition numbers associated with the Reconstruction Approach 1 increase noticeably, whereas those of the second reconstruction remain well controlled.

\begin{table}[htb]
\centering
\caption{Largest condition number of local mass matrices among all interface elements on a \xz{structured} $10\times 10$ triangular mesh.}
\label{tab:condMlocal}
\resizebox{\textwidth}{!}{%
\begin{tabular}{cccccccccc}
\toprule
\multicolumn{10}{c}{$\max\big\{\text{Cond}(M_K):~K\in\mathcal{T}_h^i\big\}$} \\
\midrule
$m$ & \multicolumn{3}{c}{$\beta^+=10$} & \multicolumn{3}{c}{$\beta^+=100$} & \multicolumn{3}{c}{$\beta^+=1000$} \\
\cmidrule(lr){2-4}\cmidrule(lr){5-7}\cmidrule(lr){8-10}
 & General & RA1 & RA2 & General & RA 1 & RA 2 & General & RA 1 & RA 2 \\
\midrule
1 & 5.42E+02 & 1        & 1 & 5.25E+04 & 1        & 1 & 5.25E+06 & 1        & 1 \\
2 & 1.19E+05 & 1        & 1 & 1.24E+07 & 1        & 1 & 1.25E+09 & 1        & 1 \\
3 & 1.47E+07 & 1        & 1 & 1.56E+09 & 1        & 1 & 1.57E+11 & 1        & 1 \\
4 & 1.06E+09 & 1        & 1 & 1.25E+11 & 1        & 1 & 1.28E+13 & 1        & 1 \\
5 & 4.18E+11 & 1        & 1 & 5.89E+13 & 1        & 1 & 6.58E+15 & 1        & 1 \\
6 & 4.13E+13 & 1        & 1 & 2.67E+15 & 1.18        & 1 & 2.60E+17 & 2.59E+01        & 1 \\
7 & 4.20E+16 & 1.27E+01        & 1 & 1.15E+19 & 1.04E+04        & 1 & 3.07E+19 & 2.06E+05        & 1 \\
8 & 6.22E+18 & 2.22E+04        & 1 & 2.56E+19 & 2.02E+07     & 1 & 3.31E+21 & 2.91E+08 & 1 \\
9 & 2.89E+21 & 1.55E+07 & 1 & 6.49E+21 & 5.23E+08 & 1 & 3.75E+22 & 1.24E+11 & 1 \\
\bottomrule
\end{tabular}%
}
\end{table}
Table~\ref{tab:condMglobal} reports the conditioning of the global mass matrices. Since the global mass matrix is block diagonal, its condition number equals the maximum condition number among all elements in the mesh. For each polynomial degree $1\le m\le 9$, we report the condition numbers on three different meshes. When \xz{$N=10$, i.e., using the structured $10\times 10$ triangular mesh}, the condition number obtained from the general construction coincides with that in Table~\ref{tab:condMlocal}, since the largest condition number occurs on the interface elements. Moreover, no apparent growth is observed under mesh refinement. This indicates that the geometric configuration of the interface elements has a much stronger impact on the conditioning than the element size itself. 

In contrast, the two reconstruction approaches yield identical condition numbers across all meshes for degrees up to $m=5$. In these cases, the condition numbers are determined entirely by non-interface elements, where an orthogonal basis \cite{2008HesthavenWarburton} is employed. For $m\ge 6$, the condition number of the Reconstruction Approach 1 is dominated by the interface elements, whereas the Reconstruction Approach 2 remains robust. This numerical experiment confirms the Remark \ref{rem: condition}.
 
\xz{Table~\ref{tab:condSglobal} reports the conditioning of the global stiffness matrices $S$ whose element components are defined in \eqref{eq: stiff local}. Although the reconstruction approaches are designed to improve the conditioning of the mass matrices, they also significantly reduce the condition number of stiffness matrices.}

\begin{table}[!htb]
\centering
\caption{Condition numbers of the global mass matrix $M$ on structured $N\times N$ triangular meshes.}
\label{tab:condMglobal}
\resizebox{\textwidth}{!}{%
\begin{tabular}{ccccccccccc}
\toprule
\multicolumn{11}{c}{Cond($M$)} \\
\midrule
$m$ & $N$ & \multicolumn{3}{c}{$\beta^+=10$} & \multicolumn{3}{c}{$\beta^+=100$} & \multicolumn{3}{c}{$\beta^+=1000$} \\
\cmidrule(lr){3-5}\cmidrule(lr){6-8}\cmidrule(lr){9-11}
 & & General & Rec 1 & Rec 2 & General & Rec 1 & Rec 2 & General & Rec 1 & Rec 2 \\
\midrule
1 & 10 & 5.4E+02 & 4  & 4  & 5.3E+04 & 4  & 4  & 5.3E+06 & 4  & 4 \\
  & 20 & 4.7E+02 & 4  & 4  & 4.6E+04 & 4  & 4  & 4.6E+06 & 4  & 4 \\
  & 40 & 3.7E+02 & 4  & 4  & 3.6E+03 & 4  & 4  & 3.6E+06 & 4  & 4 \\
\midrule
2 & 10 & 1.2E+05 & 17 & 17 & 1.2E+07 & 17 & 17 & 1.3E+09 & 17 & 17 \\
  & 20 & 3.0E+05 & 17 & 17 & 4.0E+07 & 17 & 17 & 4.1E+09 & 17 & 17 \\
  & 40 & 3.5E+05 & 17 & 17 & 4.9E+07 & 17 & 17 & 5.0E+09 & 17 & 17 \\
\midrule
3 & 10 & 1.5E+07 & 35 & 35 & 1.6E+09 & 35 & 35 & 1.6E+11 & 35 & 35 \\
  & 20 & 1.1E+08 & 35 & 35 & 1.6E+10 & 35 & 35 & 1.6E+12 & 35 & 35 \\
  & 40 & 1.4E+08 & 35 & 35 & 2.0E+10 & 35 & 35 & 2.1E+12 & 35 & 35 \\
\midrule
4 & 10 & 1.0E+09 & 46 & 46 & 1.3E+11 & 46 & 46 & 1.3E+13 & 46 & 46 \\
  & 20 & 6.5E+10 & 46 & 46 & 9.5E+12 & 46 & 46 & 9.8E+14 & 46 & 46 \\
  & 40 & 8.3E+10 & 46 & 46 & 1.2E+13 & 46 & 46 & 1.3E+15 & 46 & 46 \\
\midrule
5 & 10 & 4.2E+11 & 62 & 62 & 5.9E+13 & 62 & 62 & 6.6E+15 & 62 & 62 \\
  & 20 & 1.1E+14 & 62 & 62 & 1.7E+16 & 62 & 62 & 1.8E+18 & 62 & 62 \\
  & 40 & 1.6E+14 & 62 & 62 & 2.3E+16 & 62 & 62 & 2.4E+18 & 62 & 62 \\
\midrule
6 & 10 & 2.4E+14 & 92 & 92 & 4.5E+16 & 92 & 92 & 6.4E+18 & 3.8E+03 & 92 \\
  & 12 & 2.1E+16 & 92 & 92 & 3.5E+18 & 2.2E+02 & 92 & 4.6E+20 & 3.2E+04 & 92 \\
  & 14 & 4.8E+16 & 92 & 92 & 6.7E+18 & 6.0E+02 & 92 & 8.4E+20 & 6.3E+04 & 92 \\
\midrule
7 & 10 & 4.2E+16 & 125 & 125 & 1.2E+19 & 1.0E+04 & 125 & 3.1E+19 & 2.1E+05 & 125 \\
  & 12 & 5.4E+18 & 1.1E+03 & 125 & 7.5E+20 & 3.8E+05 & 125 & 1.3E+23 & 2.3E+07 & 125 \\
  & 14 & 3.7E+19 & 1.0E+04 & 125 & 4.8E+21 & 3.7E+06 & 125 & 1.7E+24 & 2.5E+08 & 125 \\
\midrule
8 & 10 & 6.2E+18 & 2.2E+04 & 193 & 2.6E+19 & 2.0E+06 & 193 & 3.3E+21 & 2.9E+08 & 193 \\
  & 12 & 1.4E+22 & 1.8E+07 & 193 & 2.3E+25 & 1.5E+09 & 193 & 1.1E+25 & 1.2E+11 & 193 \\
  & 14 & 2.7E+22 & 1.3E+07 & 193 & 2.6E+24 & 2.2E+09 & 193 & 2.9E+26 & 4.7E+11 & 193 \\
\midrule
9 & 10 & 2.9E+21 & 1.6E+07 & 285 & 6.5E+21 & 5.2E+08 & 285 & 3.8E+22 & 1.2E+11 & 285 \\
  & 12 & 9.7E+23 & 2.3E+09 & 285 & 5.5E+25 & 3.4E+11 & 285 & 1.9E+26 & 2.1E+13 & 285 \\
  & 14 & 3.9E+24 & 5.3E+09 & 285 & 6.0E+26 & 1.4E+12 & 285 & 5.4E+29 & 1.0E+14 & 285 \\
\bottomrule
\end{tabular}%
}
\end{table}

\begin{table}[!htb]
\centering
\caption{Condition numbers of the global stiffness matrices on structured $N\times N$ triangular meshes.}
\label{tab:condSglobal}
\resizebox{\textwidth}{!}{%
\begin{tabular}{ccccccccccc}
\toprule
\multicolumn{11}{c}{Cond($S$)} \\
\midrule
$m$ & $N$ & \multicolumn{3}{c}{$\beta^+=10$} & \multicolumn{3}{c}{$\beta^+=100$} & \multicolumn{3}{c}{$\beta^+=1000$} \\
\cmidrule(lr){3-5}\cmidrule(lr){6-8}\cmidrule(lr){9-11}
 & & General & Rec 1 & Rec 2 & General & Rec 1 & Rec 2 & General & Rec 1 & Rec 2 \\
\midrule
1 & 10 & 3.7E+01 & 1.1E+01  & 1.1E+01  & 2.0E+03 & 1.1E+02  & 1.1E+02  & 4.0E+04 & 1.1E+03  & 1.1E+03 \\
   & 20 & 7.8E+01 & 1.1E+01  & 1.1E+01  & 5.0E+03 & 1.0E+02  & 1.0E+02  & 1.1E+05 & 1.0E+03  & 1.0E+03 \\
   & 40 & 7.1E+01 & 1.2E+01  & 1.2E+01  & 2.3E+03 & 1.1E+02  & 1.1E+02  & 3.0E+04 & 1.1E+03  & 1.1E+03 \\
\midrule
2 & 10 & 4.9E+03 & 6.0E+01  & 6.0E+01  & 4.4E+05 & 5.0E+02  & 5.0E+02  & 3.5E+07 & 4.9E+03  & 4.9E+03 \\
   & 20 & 1.1E+04 & 8.1E+01  & 8.1E+01  & 1.3E+06 & 7.9E+02  & 7.9E+02  & 1.2E+08 & 7.8E+03  & 7.8E+03 \\
   & 40 & 1.2E+04 & 9.8E+01  & 9.8E+01  & 1.5E+06 & 9.7E+02  & 9.7E+02  & 1.0E+08 & 9.7E+03  & 9.7E+03 \\
\midrule
3 & 10 & 4.6E+05 & 2.8E+02  & 2.8E+02  & 4.3E+07 & 2.7E+03  & 2.7E+03  & 3.2E+09 & 2.7E+04  & 2.7E+04 \\
   & 20 & 7.1E+06 & 2.8E+02  & 2.8E+02  & 9.5E+08 & 2.7E+03  & 2.7E+03  & 7.6E+10 & 2.7E+04  & 2.7E+04 \\
   & 40 & 8.5E+06 & 1.6E+02  & 1.6E+02  & 9.4E+08 & 1.5E+03  & 1.5E+03  & 7.9E+10 & 1.5E+04  & 1.5E+04 \\
\midrule
4 & 10 & 7.4E+07 & 4.8E+02  & 4.8E+02  & 8.1E+09 & 4.4E+03  & 4.4E+03  & 6.5E+11 & 4.4E+04  & 4.4E+04 \\
   & 20 & 4.3E+09 & 4.7E+02  & 4.7E+02  & 6.0E+11 & 4.6E+03  & 4.6E+03  & 4.8E+13 & 4.6E+04  & 4.6E+04 \\
   & 40 & 4.7E+09 & 6.3E+02  & 6.3E+02  & 5.2E+11 & 6.1E+03  & 6.1E+03  & 2.2E+13 & 6.1E+04  & 6.1E+04 \\
\midrule
5 & 10 & 1.7E+10 & 7.8E+02  & 7.8E+02  & 2.0E+12 & 6.4E+03  & 6.4E+03  & 1.8E+14 & 6.3E+04  & 6.3E+04 \\
   & 20 & 6.8E+12 & 9.4E+02  & 9.4E+02  & 9.2E+13 & 8.7E+03  & 8.7E+03  & 8.4E+13 & 8.7E+04  & 8.7E+04 \\
   & 40 & 4.3E+12 & 1.4E+03  & 1.4E+03  & 1.7E+14 & 1.3E+04  & 1.3E+04  & 6.1E+13 & 1.3E+05  & 1.3E+05 \\
\midrule
6 & 10 & 8.2E+12 & 1.7E+03  & 1.7E+03  & 1.0E+14 & 1.4E+04  & 1.4E+04  & 2.5E+13 & 2.7E+05  & 1.4E+05 \\
   & 12 & 2.1E+14 & 1.6E+03  & 1.5E+03  & 2.3E+14 & 2.2E+04  & 1.3E+04  & 2.8E+14 & 5.9E+05  & 1.3E+05 \\
   & 14 & 2.0E+14 & 2.4E+03  & 2.4E+03  & 1.9E+14 & 1.4E+04  & 2.4E+04  & 1.0E+14 & 8.3E+05  & 2.4E+05 \\
\midrule

7 & 10 & 1.9E+14 & 2.5E+03  & 2.6E+03  & 9.2E+13 & 1.3E+05  & 2.5E+04  & 9.3E+13 & 5.8E+07  & 2.6E+05 \\
   & 12 & 8.5E+13 & 2.6E+03  & 2.8E+03  & 1.3E+14 & 2.2E+06  & 2.9E+04  & 1.8E+14 & 5.3E+07  & 2.9E+05 \\
   & 14 & 1.3E+14 & 1.3E+04  & 4.7E+03  & 1.3E+14 & 6.3E+06  & 4.5E+04  & 2.3E+14 & 5.6E+08  & 4.4E+05 \\
\midrule

8 & 10 & 8.4E+13 & 6.5E+04  & 4.8E+03  & 8.8E+13 & 1.8E+07  & 4.8E+04  & 1.1E+14 & 3.2E+10  & 4.8E+05 \\
   & 12 & 1.2E+14 & 2.1E+07  & 5.8E+03  & 1.6E+14 & 6.8E+08  & 6.2E+04  & 1.7E+14 & 5.3E+09  & 6.2E+05 \\
   & 14 & 1.3E+14 & 1.8E+07  & 5.4E+03  & 1.5E+14 & 8.1E+08  & 5.5E+04  & 1.7E+14 & 3.5E+10  & 5.5E+05 \\
\midrule

9 & 10 & 9.6E+13 & 3.5E+07  & 6.6E+03  & 1.3E+14 & 1.0E+10  & 6.1E+04  & 1.3E+14 & 1.2E+12  & 6.0E+05 \\
   & 12 & 1.1E+14 & 9.0E+08  & 1.0E+04  & 1.4E+14 & 3.8E+10  & 1.0E+05  & 1.3E+14 & 6.2E+11  & 1.0E+06 \\
   & 14 & 1.3E+14 & 3.7E+09  & 6.9E+03  & 1.2E+14 & 6.2E+10  & 6.7E+04  & 1.4E+14 & 8.5E+11  & 6.7E+05 \\
\bottomrule
\end{tabular}%
}
\end{table}

\subsection{Numerical Experiments on $L^2$ Projection}	
\begin{figure}[h]
\includegraphics[width = .49\textwidth]{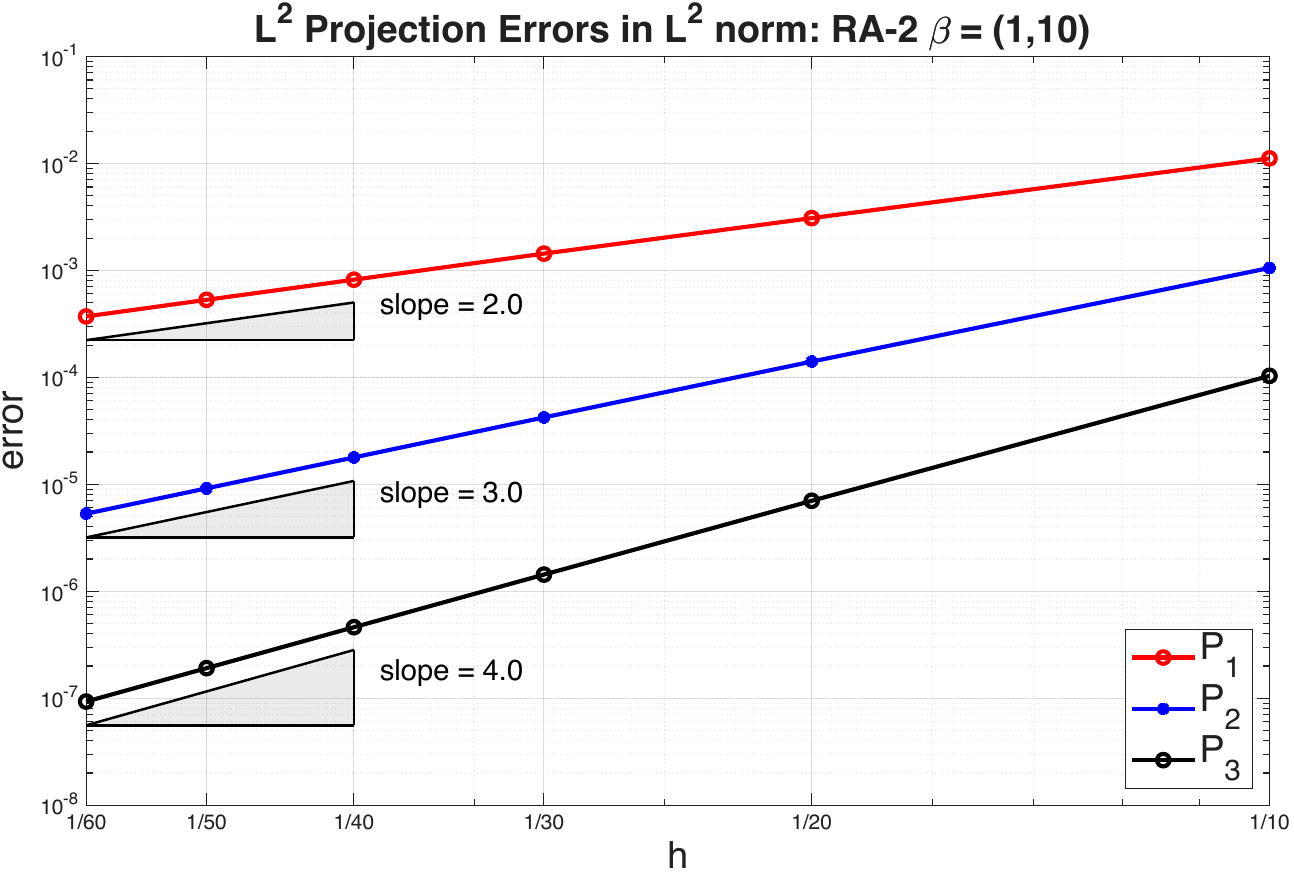}
\includegraphics[width = .49\textwidth]{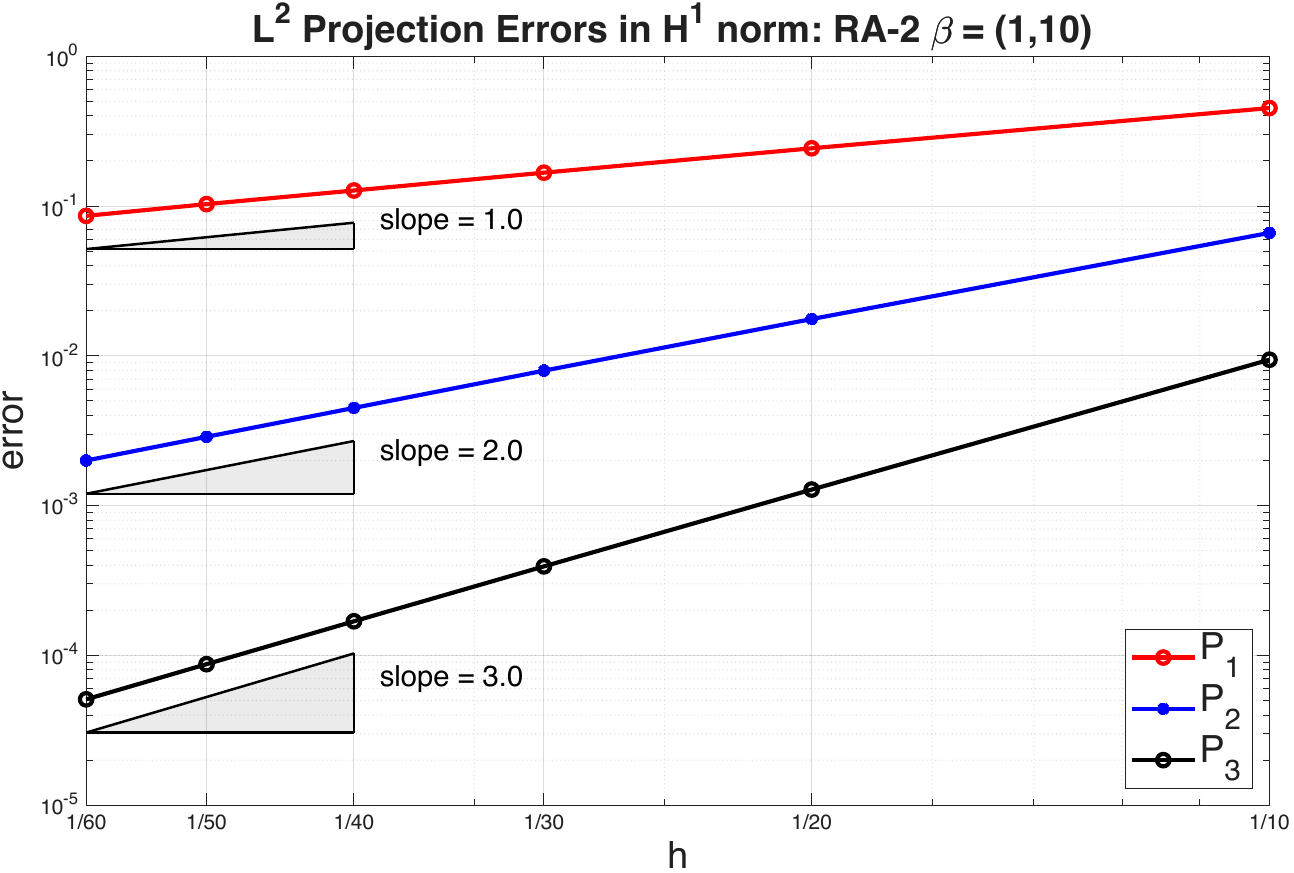}
\caption{Convergence of $L^2$ projection in $L^2$ and $H^1$ norms using Reconstruction Approach 2 with degree $1\le m\le 3$.}
\label{fig: proj Re2 p1-p3 1-10}
\end{figure}

\begin{figure}[h]
\includegraphics[width = .49\textwidth]{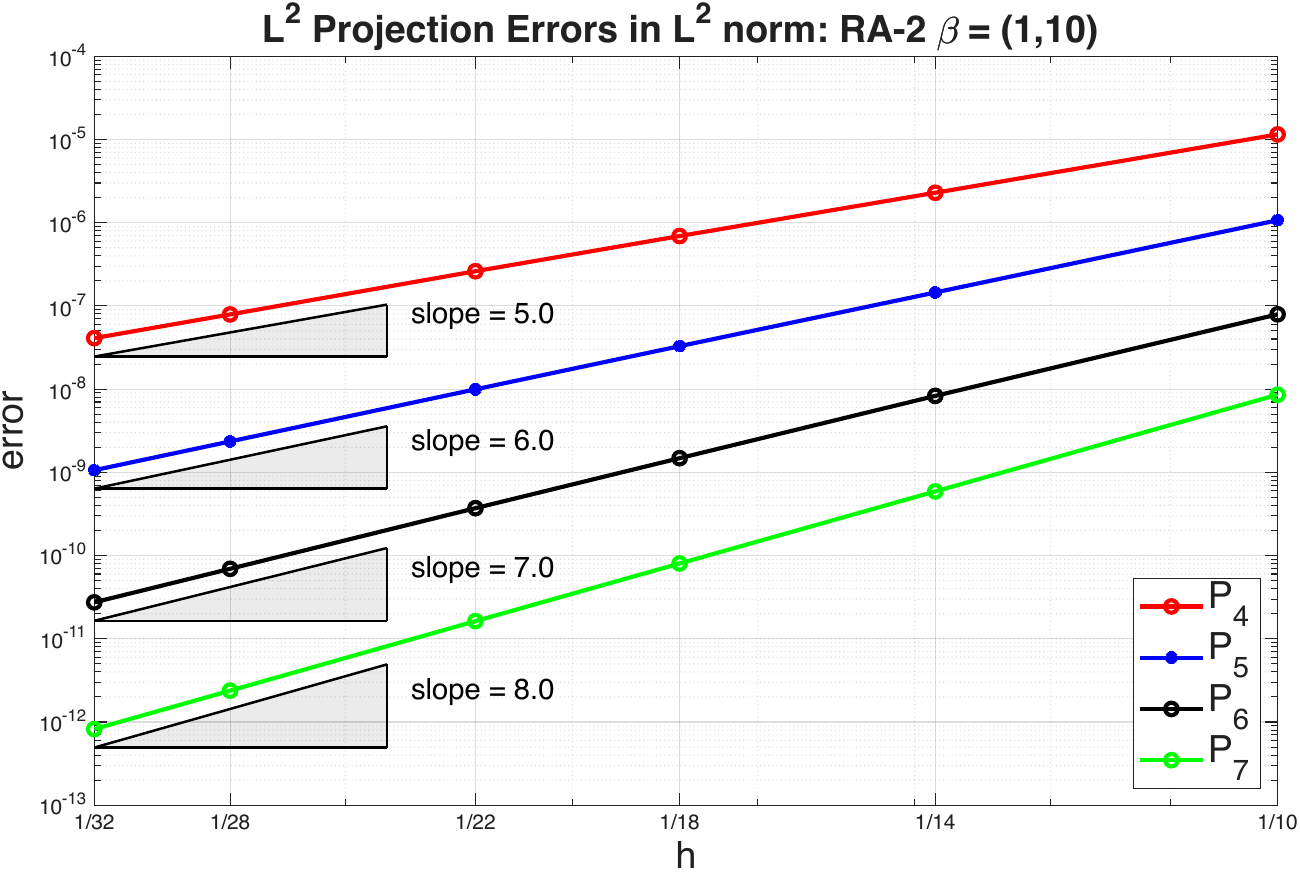}
\includegraphics[width = .49\textwidth]{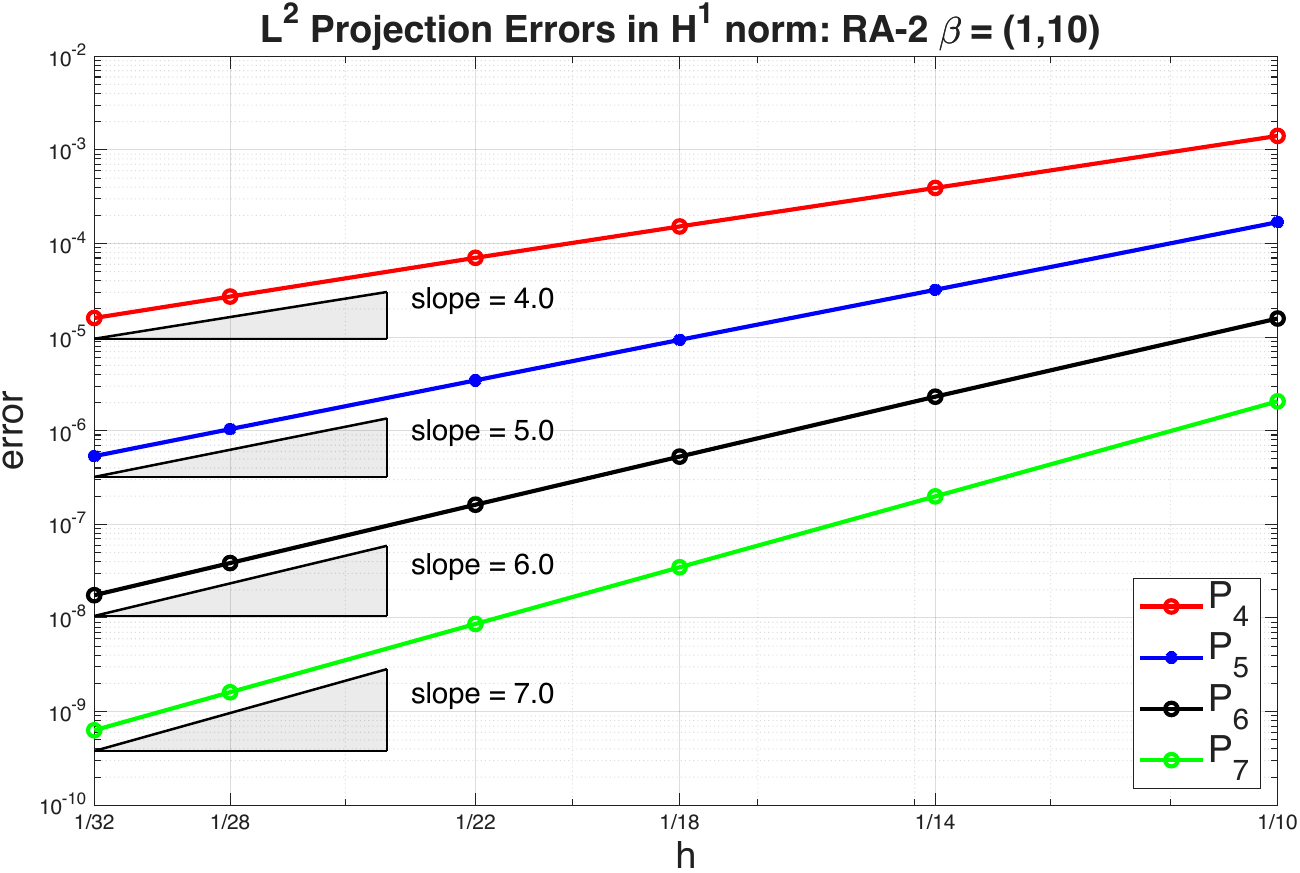}
\caption{Convergence of $L^2$ projection in $L^2$ and $H^1$ norms using Reconstruction Approach 2 with degree $4\le m\le 7$.}
\label{fig: proj Re2 p4-p7 1-10}
\end{figure}

In this experiment, we examine the approximation properties of the GC-IFE spaces through the convergence of the $L^2$
projection $\mathcal{P}_hu$. All three types of IFE constructions are tested. For polynomial degrees $m\le 3$, the convergence results obtained using Reconstruction Approach 2 are presented in Figure~\ref{fig: proj Re2 p1-p3 1-10}, in both the 
$L^2$ and semi-$H^1$-norms. For higher polynomial degrees $4\le m\le 7$, we use coarser meshes  to minimize the impact from round-off errors. The convergence results using Reconstruction Approach 2 are presented in Figure~\ref{fig: proj Re2 p4-p7 1-10}. In both cases, the behaviors of the general construction and Reconstruction Approach 1 are very similar, and therefore their results are omitted for brevity. The observed convergence rates indicate that
 \[
\|u-\mathcal{P}_hu\|_{L^2(\Omega)} \approx \mathcal{O}(h^{m+1}),\qquad
|u-\mathcal{P}_hu|_{H^1(\Omega)} \approx \mathcal{O}(h^{m}).
\]

\begin{figure}[h]

\includegraphics[width = .49\textwidth]{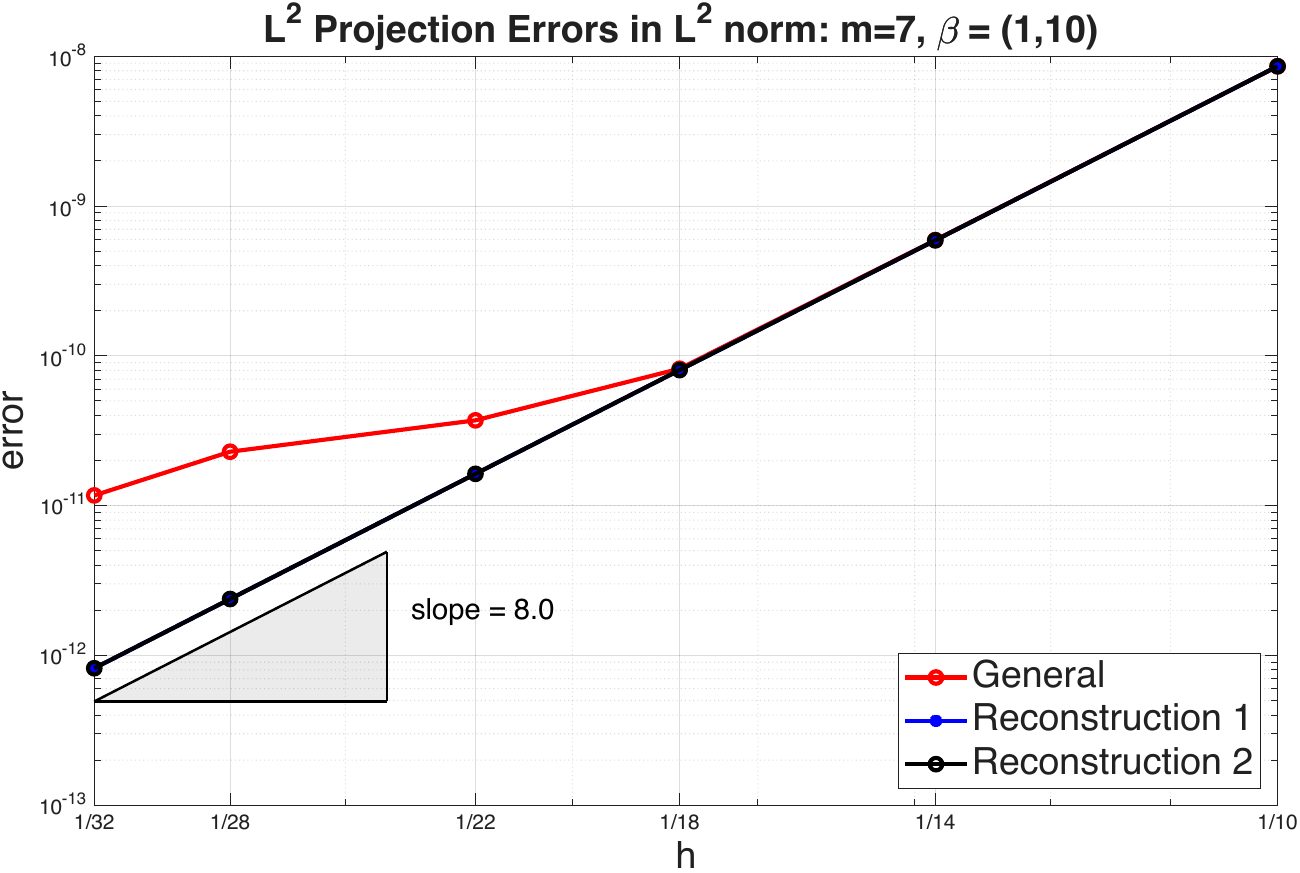}
\includegraphics[width = .49\textwidth]{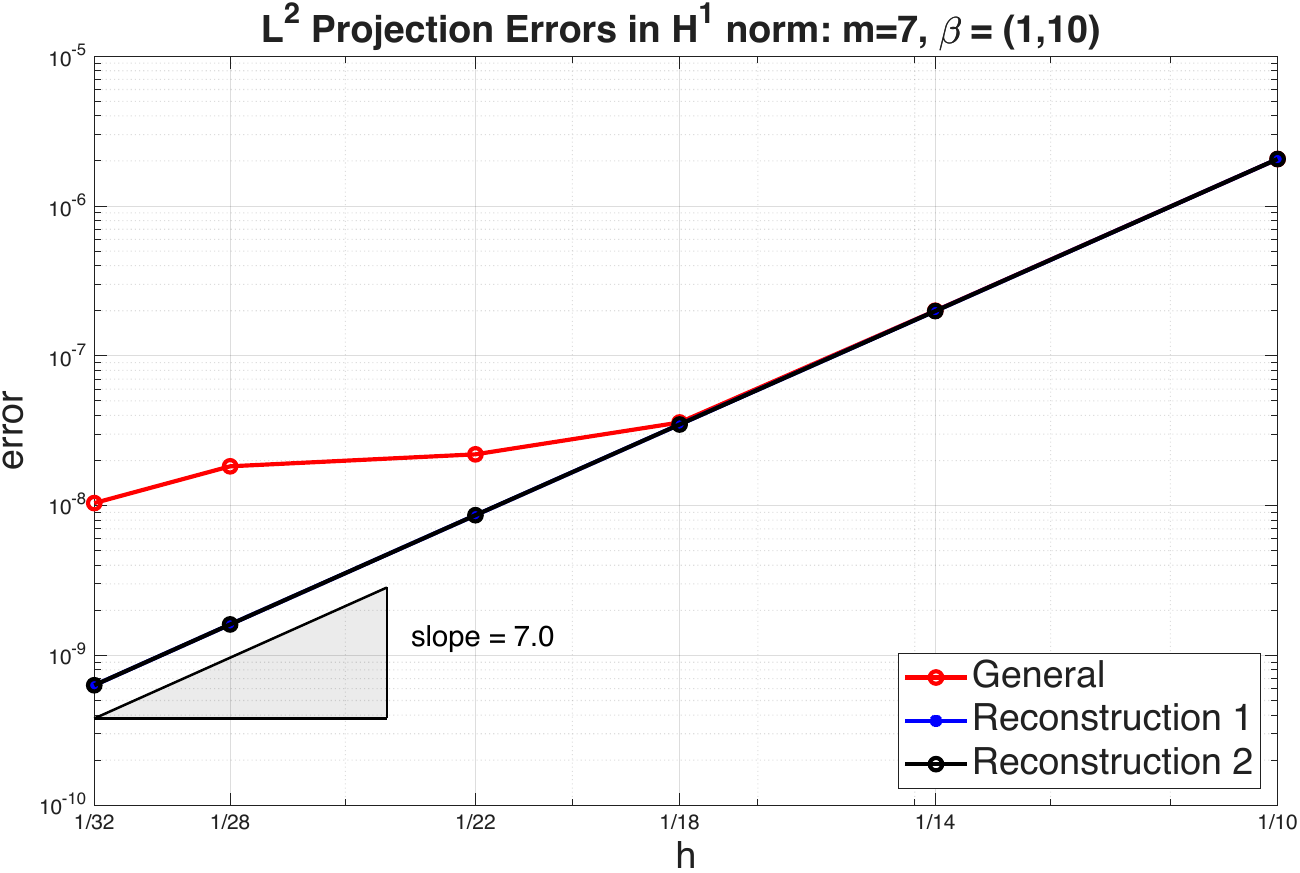}
\caption{Convergences of $L^2$ projection error in $L^2$ and $H^1$ norms using different bases with degree $m=7$.}
\label{fig: proj p7 1-10}
\end{figure}

As discussed in Section~\ref{sec: condition}, for higher-order approximations (e.g. $m=7$), the conditioning of the mass matrices differs significantly among the three bases. As shown in Figure~\ref{fig: proj p7 1-10}, the accuracy of the approximation based on the general construction deteriorates significantly, whereas both reconstructed approaches continue to exhibit convergence.

	\medskip
\subsection{IFE-DG Solution for Elliptic Interface Problems}	
\begin{figure}[h]
\includegraphics[width = .49\textwidth]{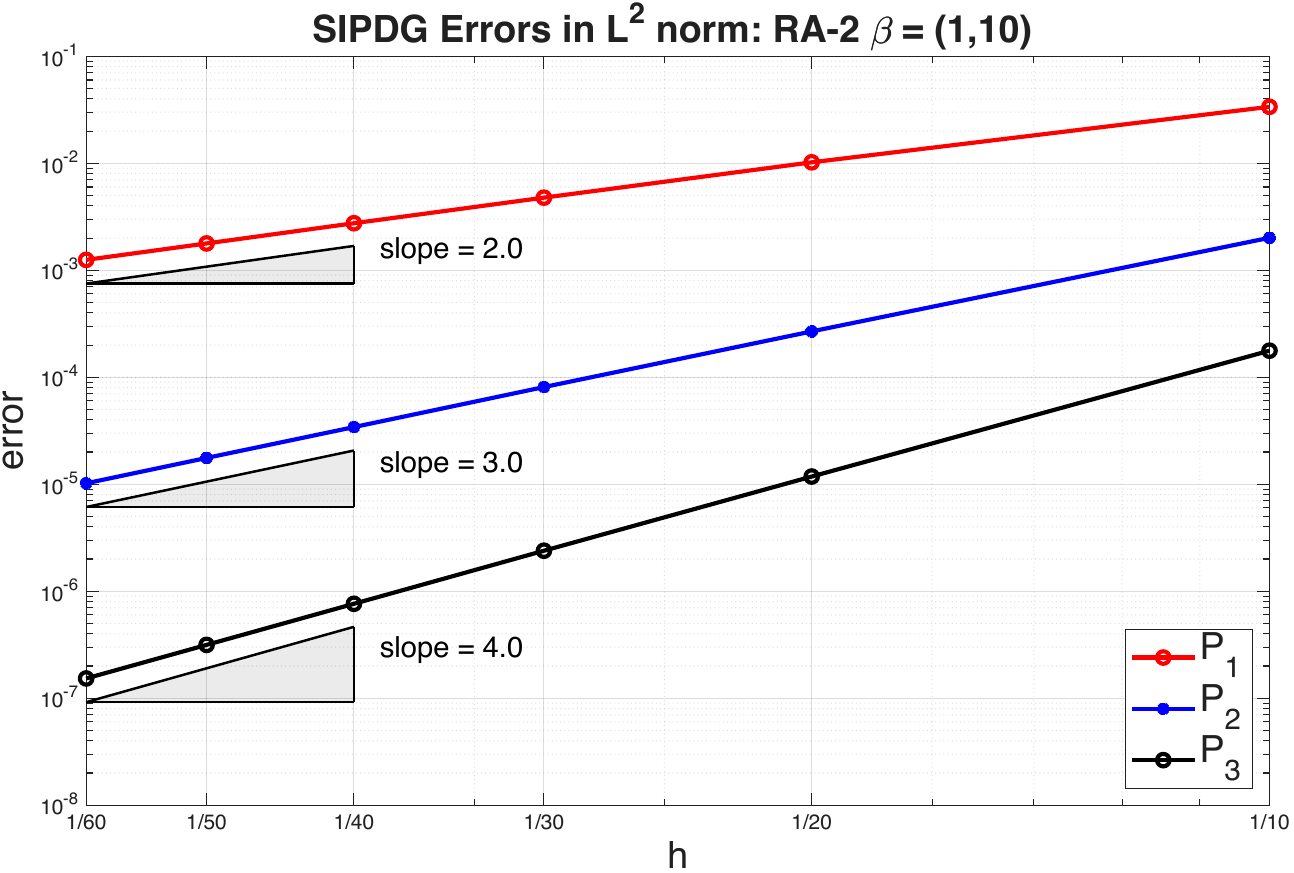}
\includegraphics[width = .49\textwidth]{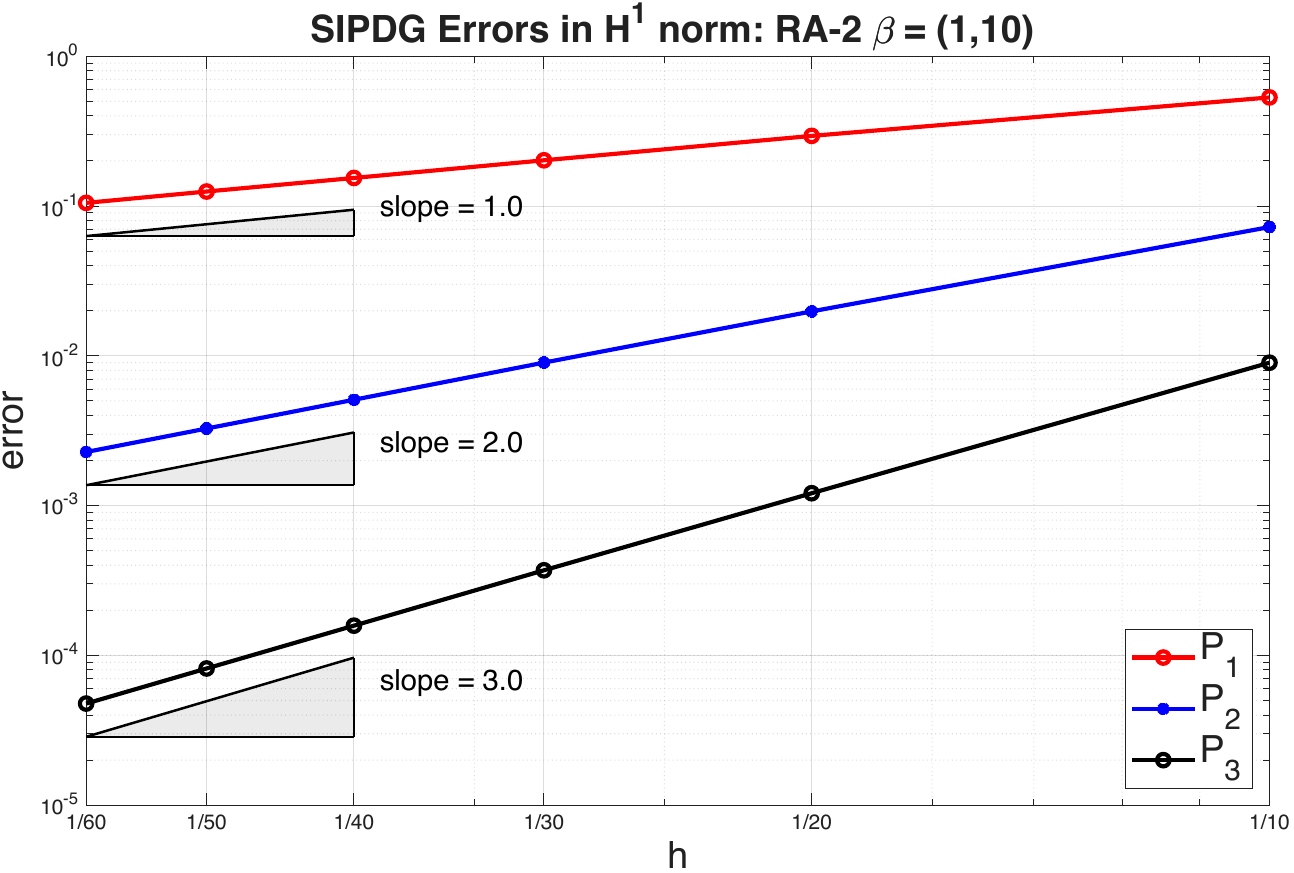}
\caption{Convergence of SIPDG solution in $L^2$ and $H^1$ norms using Reconstruction Approach 2 with degrees up to 3.}
\label{fig: sol RA2 p1-p3 1-10}
\end{figure}

\begin{figure}[h]
\includegraphics[width = .49\textwidth]{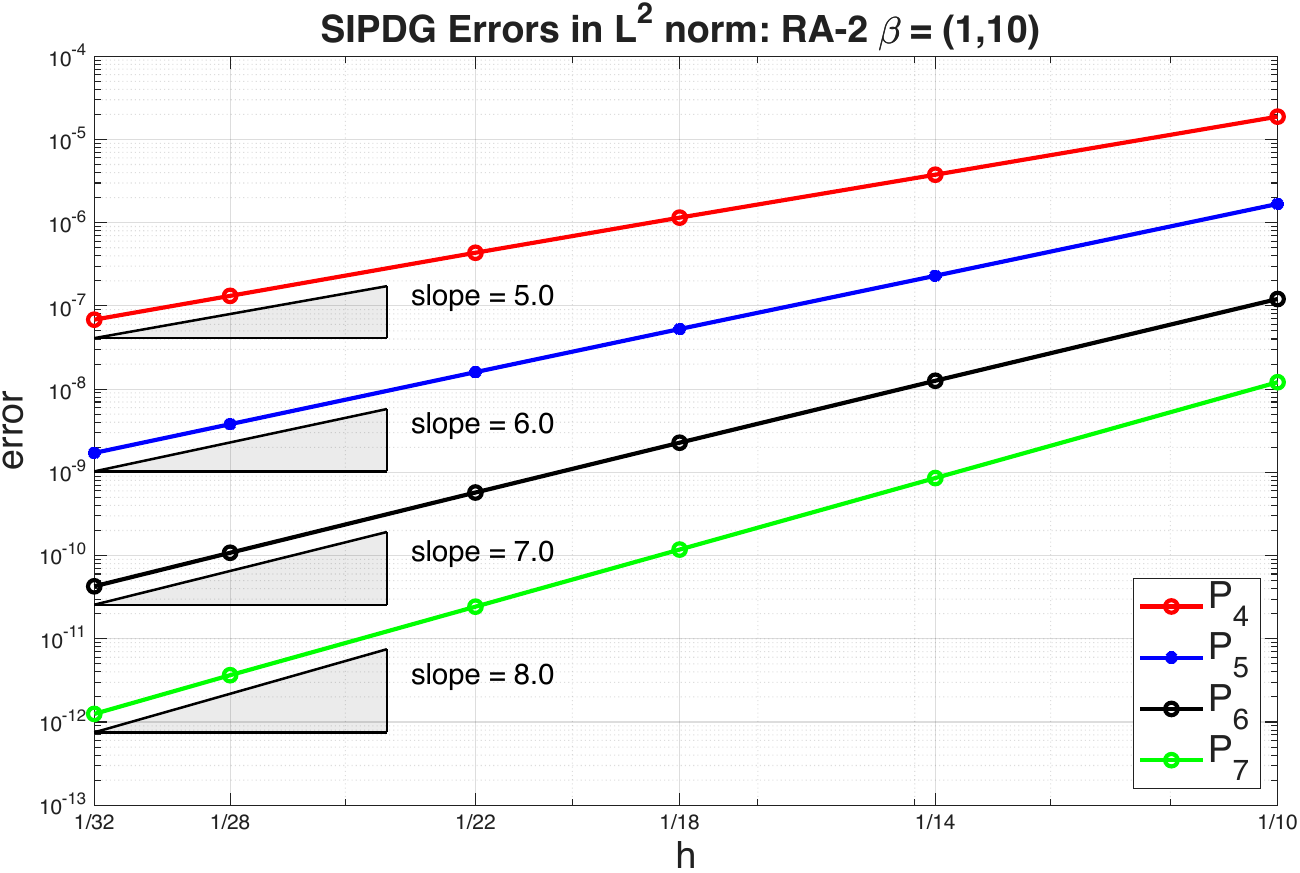}
\includegraphics[width = .49\textwidth]{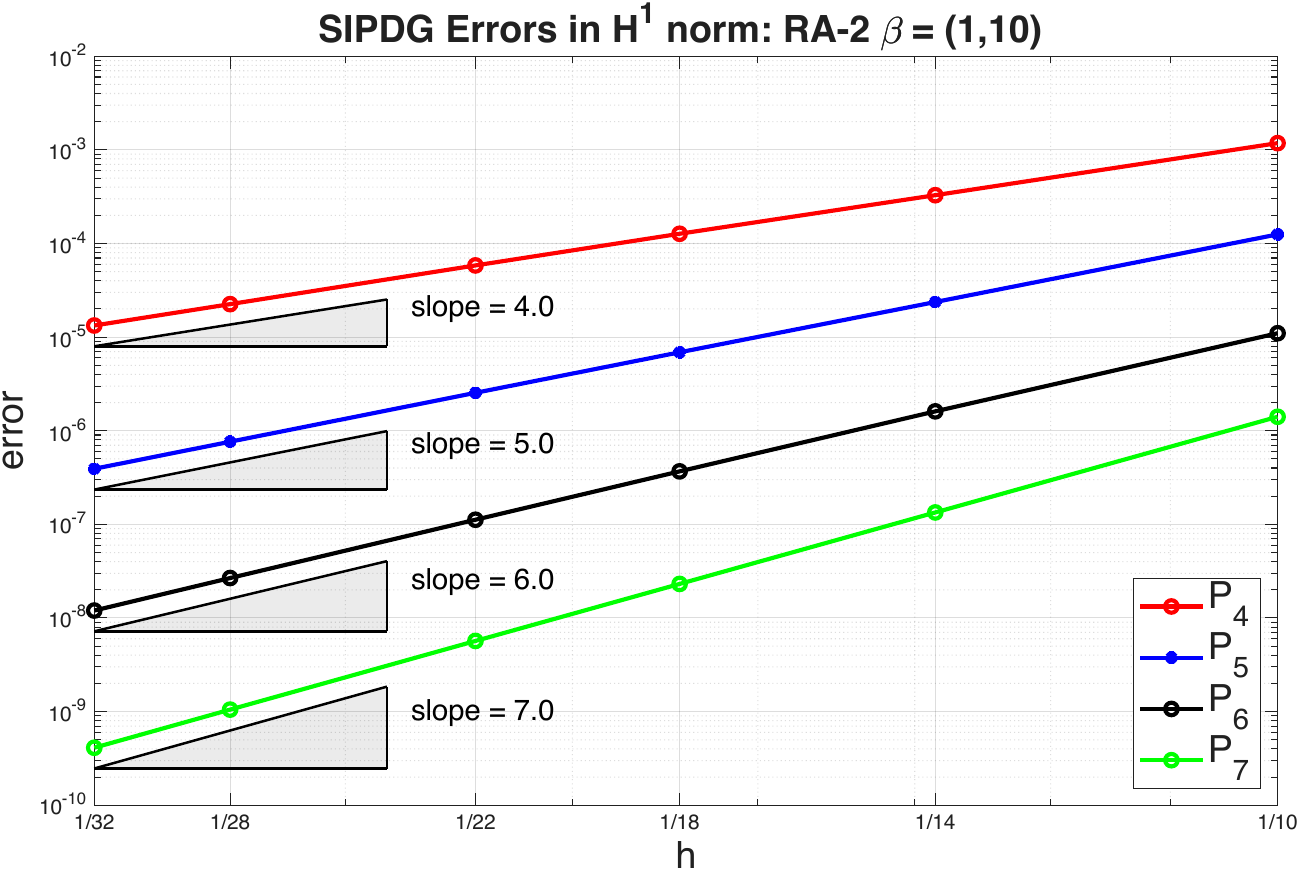}
\caption{Convergence of SIPDG solution in $L^2$ and $H^1$ norms using Reconstruction Approach 2 with degrees $4\le m\le 7$.}
\label{fig: sol RA2 p4-p7 1-10}
\end{figure}

	\smallskip
	Next, we apply our GC-IFE spaces in the symmetric interior penalty discontinuous Galerkin (SIPDG) method for solving the elliptic interface problems \eqref{eq:FrenentIFE_ellipPDE_tri}-\eqref{eq: jumpEx}.
    Find $u_h\in {\mathcal{V}}^m_{\beta}({\mathcal{T}_h})$ such that
    \begin{equation*}
        \begin{split}
            \sum_{K\in\mathcal{T}_h} \int_K \beta \nabla u_h\cdot \nabla v_h d\mathbf{x} - \sum_{e\in\mathcal{E}^i_h} \int_e \{ \beta \nabla u_h\cdot \mathbf{ n} \} \jump{v_h} ds  & \\
         - \sum_{e\in\mathcal{E}^i_h} \int_e \{ \beta \nabla v_h \cdot \mathbf{ n} \} \jump{u_h} ds + \sum_{e\in\mathcal{E}^i_h} \frac{\sigma_e }{|e|} \int_e  \jump{u_h} \, \jump{v_h} ds &= \int_{\Omega}fv_h d\mathbf{x},~ \forall v_h\in {\mathcal{V}}^m_{\beta}({\mathcal{T}_h}).
        \end{split}
    \end{equation*}
    Here, the penalty parameter $\sigma_e$ is set to be:
	\[
	\sigma_e = \sigma_0 \frac{\beta_{\max}^2}{\beta_{\min}}, \quad \sigma_0 = \xz{m^2},
	\]
     where $\beta_{\max} = \max\{\beta^+,\beta^-\}$ and $\beta_{\min} = \min\{\beta^+,\beta^-\}$. 
In Figures \ref{fig: sol RA2 p1-p3 1-10} and \ref{fig: sol RA2 p4-p7 1-10}, we present the convergence of SIPDG solutions using IFE bases from Reconstruction Approach 2. The performance of IFE basis from Reconstruction Approach 1 is very similar and is therefore omitted. These numerical results demonstrate that the DG solutions constructed with all three bases achieve the expected convergence orders in both the $L^2$ and $H^1$ norms. 

\begin{figure}[h]
\includegraphics[width = .49\textwidth]{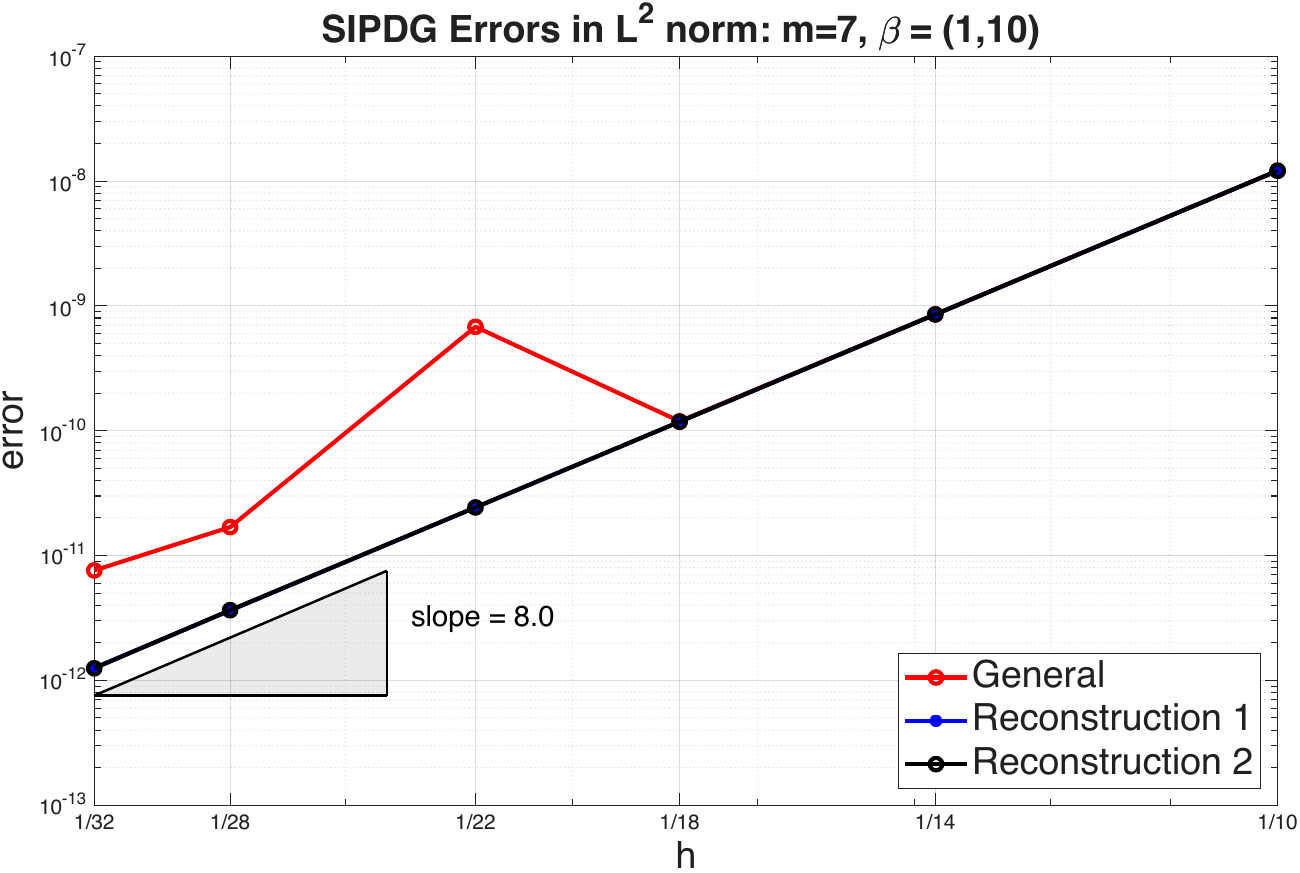}
\includegraphics[width = .49\textwidth]{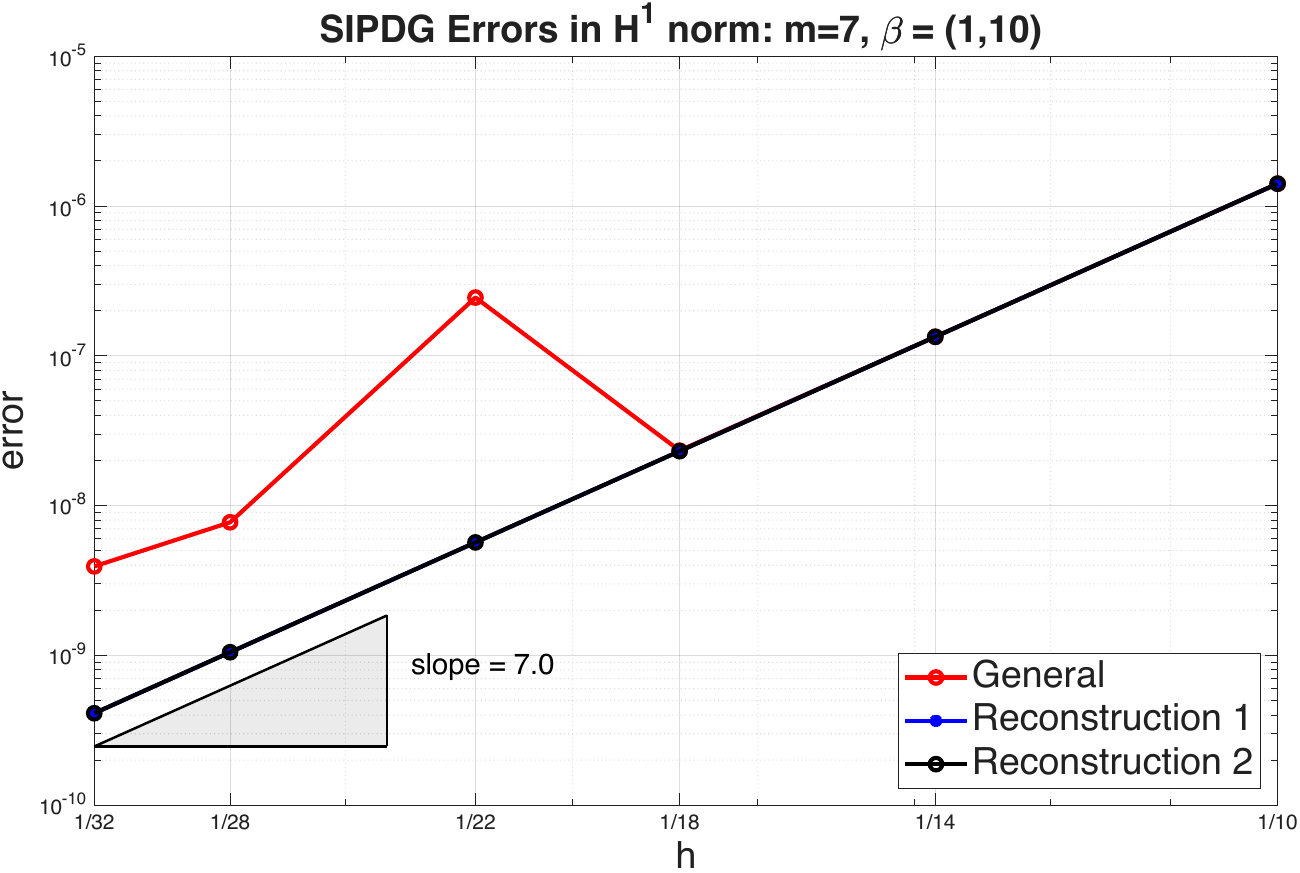}
\caption{Convergence of DG-IFE errors in $L^2$ and $H^1$ norm using all three bases with degree $m=7$.}
\label{fig: sol 1-10 7}
\end{figure}

The computational results using general construction are also similar for low-order polynomials ($m\le 6$). For higher-degree approximation (e.g. $m=7$), we report the convergence of the DG solutions generated from all three bases in Figure \ref{fig: sol 1-10 7}. Again, the reconstructed bases exhibit more robust performance compared with the general construction.

	\medskip

\begin{figure}[!htb]
\includegraphics[width = .49\textwidth]{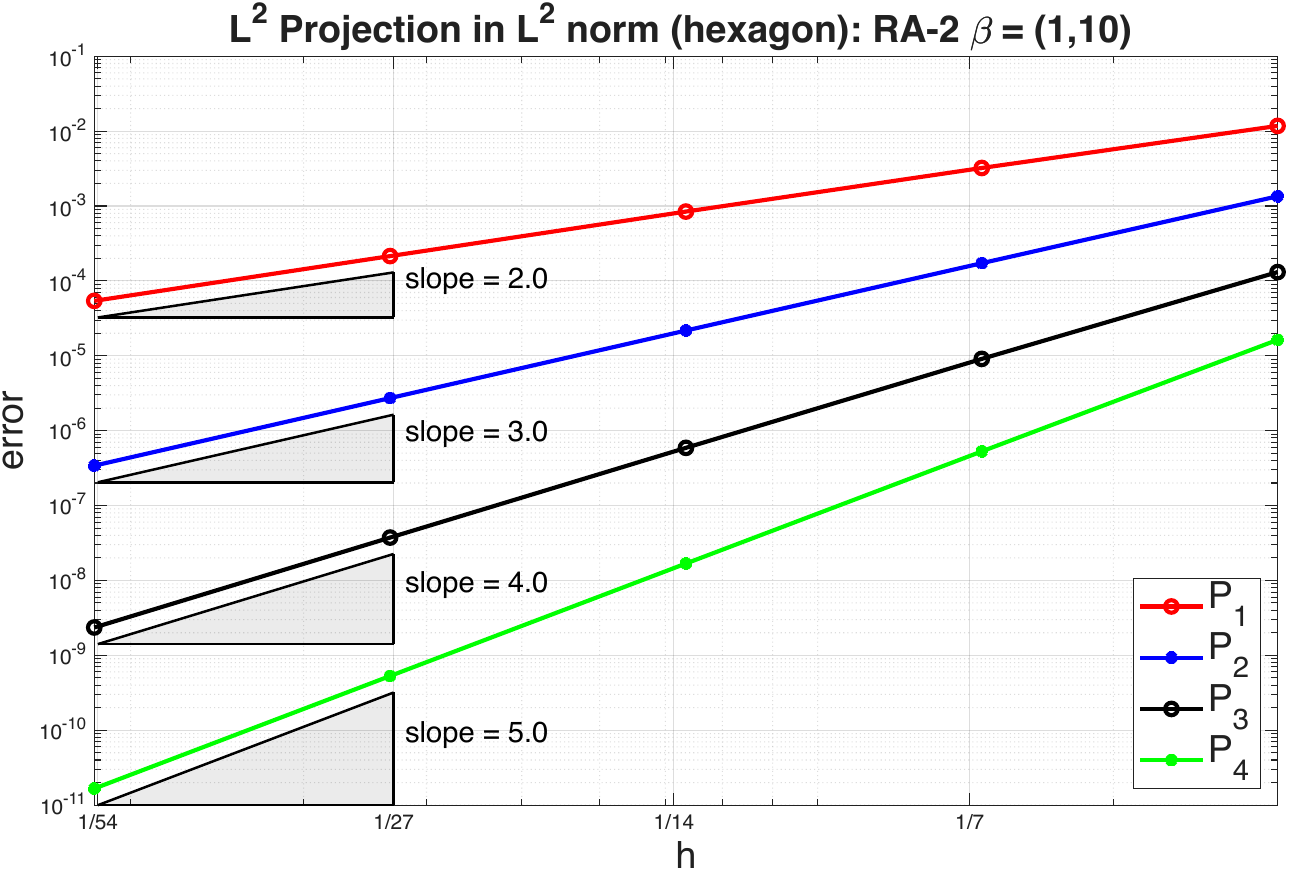}
\includegraphics[width = .49\textwidth]{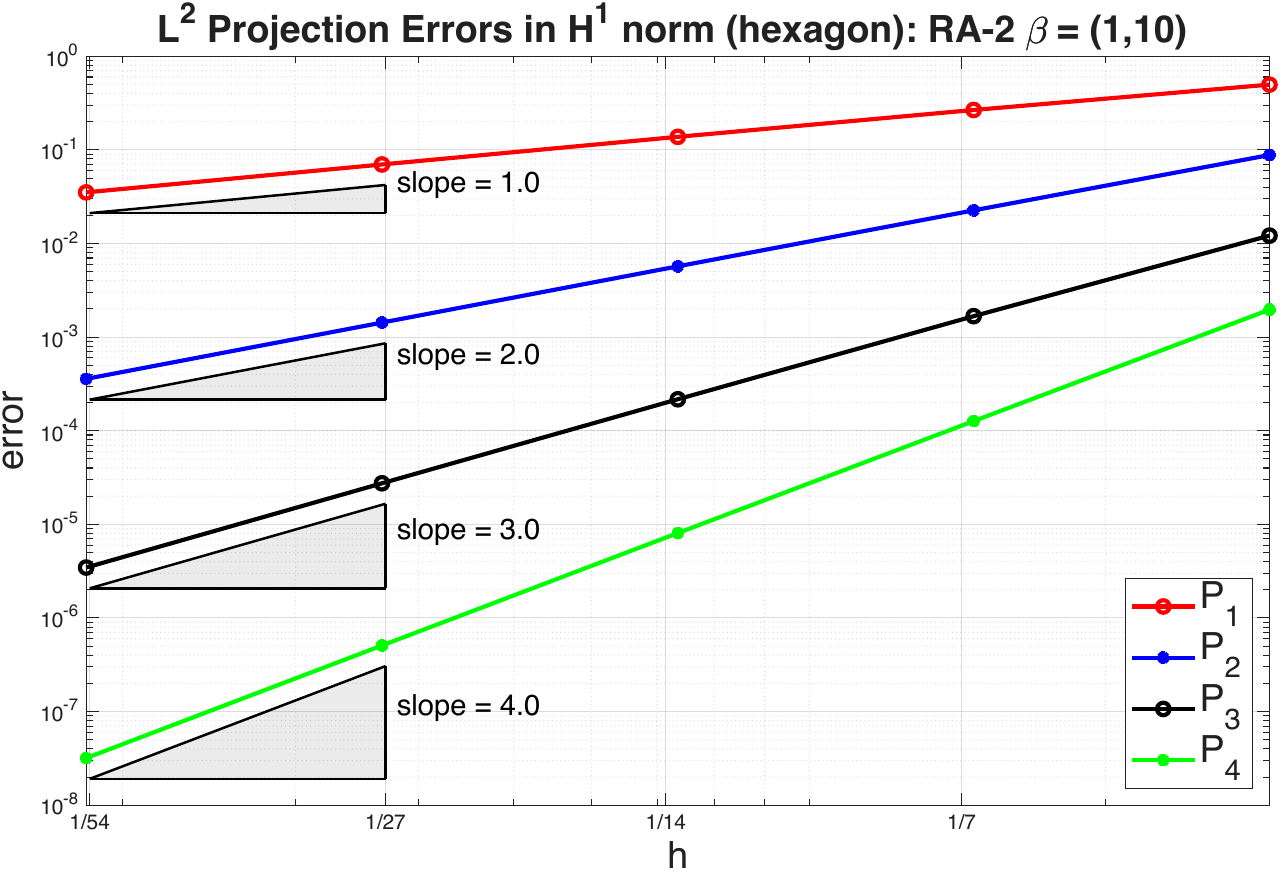}
\caption{$L^2$ Projection errors in $L^2$ and $H^1$ norms on hexagonal domain using Reconstruction Approach 2 with degrees up to $4$.}
\label{fig: proj hex 1-10}
\end{figure}

\begin{figure}[!htb]
\includegraphics[width = .49\textwidth]{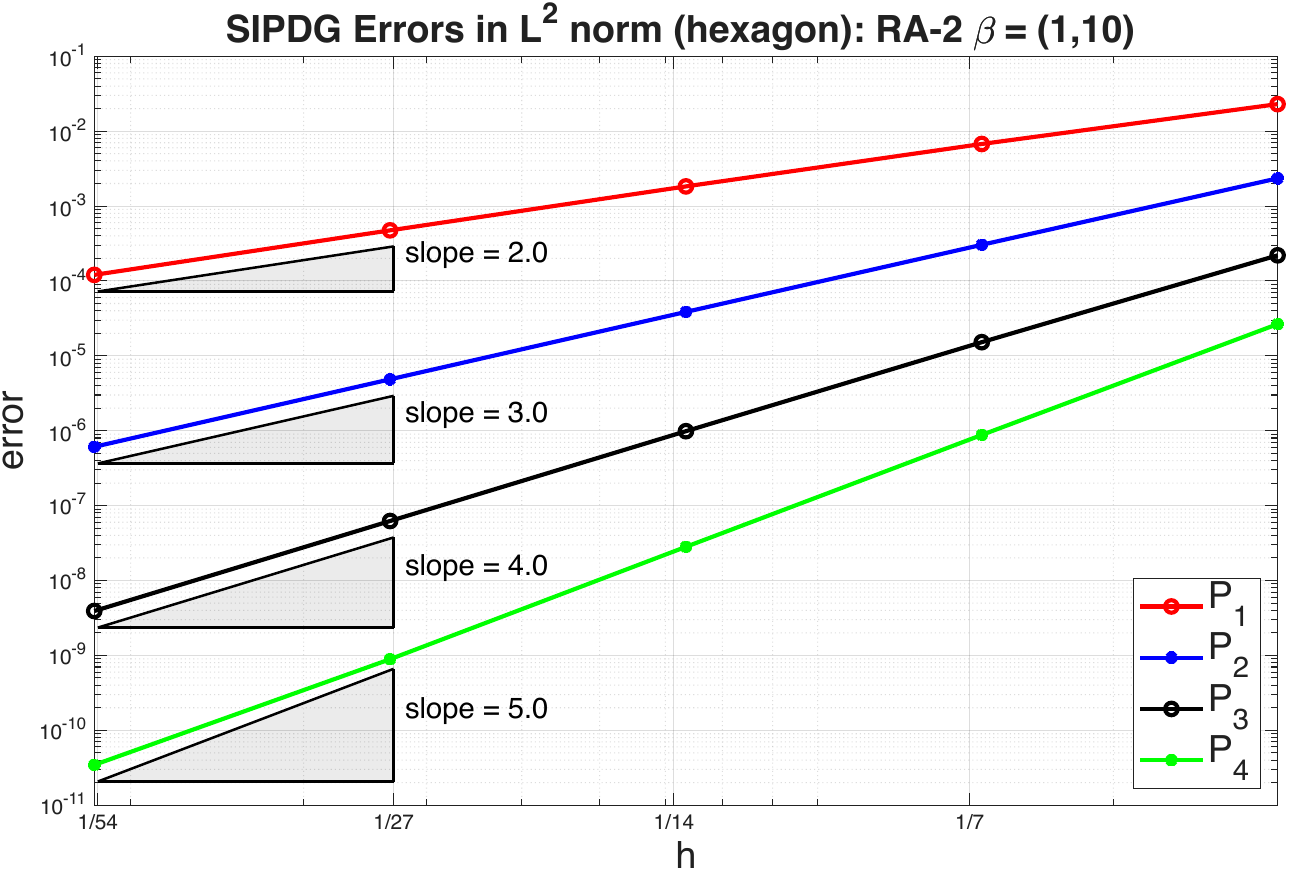}
\includegraphics[width = .49\textwidth]{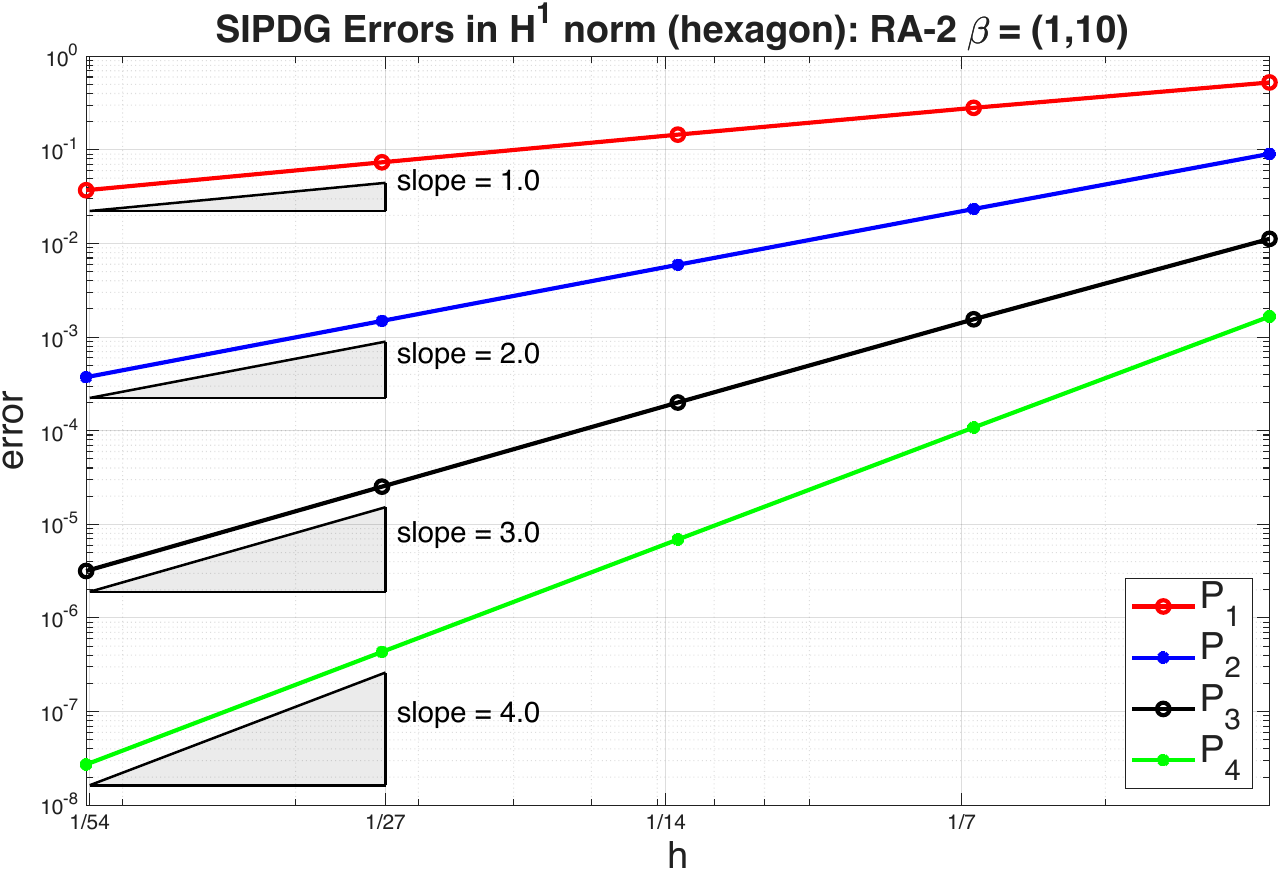}
\caption{SIPDG errors in $L^2$ and $H^1$ norms on hexagonal domain using Reconstruction Approach 2 with degrees up to $4$.}
\label{fig: solution hex 1-10}
\end{figure}

\xz{
\subsection{Non-rectangular Domain}	
In this example, we consider an interface problem defined on a non-rectangular domain. Since rectangular meshes are not applicable in this setting, we use triangular meshes to partition the domain.  The computational domain is a regular hexagon with vertices 
\[
A_k = \left(\cos\frac{(2k+1)\pi}{6},\ \sin\frac{(2k+1)\pi}{6}\right), \qquad 0 \le k \le 5.
\] 
The interface $\Gamma$ is taken to be the same circle as in previous examples, and the exact solution is defined as in \eqref{eq: true solution}. The triangular meshes are generated with MATLAB PDE Toolbox independent of the interface, as illustrated in the right plot of Figure \ref{fig: mesh}. 

The convergence of the $L^2$ projection using Reconstruction Approach 2 is reported in Figure \ref{fig: proj hex 1-10}. Optimal convergence rates are observed in both the $L^2$ norm and the semi-$H^1$ norm. The convergence of the SIPDG solution error is reported in Figure \ref{fig: solution hex 1-10}, where optimal decay rates are again achieved in both norms. 
The convergence behavior for higher-contrast coefficients is qualitatively similar and is therefore omitted. Likewise, for 
higher-degree approximations ($m \ge 5$), the results are analogous to those in 
Figure~\ref{fig: proj Re2 p4-p7 1-10} and Figure~\ref{fig: sol 1-10 7} and are omitted for brevity.
}
	
\section{Conclusion}\label{sec: conclusion}
We have developed geometry-conforming immersed finite element spaces on triangular meshes for elliptic interface problems using a Frenet–Serret mapping that straightens curved interfaces and enforces the jump conditions exactly. Extending earlier work on rectangular meshes, we introduced three high-order Frenet-IFE constructions, including reconstruction approaches that significantly improve mass matrix conditioning. Numerical experiments confirm the optimal approximation properties and robustness of the proposed spaces. When combined with interior penalty discontinuous Galerkin methods, the resulting schemes achieve optimal convergence rates in both the 
$H^1$- and $L^2$- norms. These results provide a foundation for future work on extending the Frenet mapping to handle boundary value problems with curved boundaries and high-order fitted-mesh methods for interface problems.

\section*{Acknowledgments}
Xu Zhang is partially supported by National Science Foundation grant DMS-2110833.
\bibliographystyle{amsplain}
\bibliography{XZbib}

@article{2017GuzmanSanchezSarkis,
	author = {Johnny Guzm\'an and Manuel A. S\'anchez and Marcus Sarkis},
	fjournal = {Journal of Scientific Computing},
	journal = {J. Sci. Comput.},
	number = {1},
	pages = {330--365},
	title = {A finite element method for high-contrast interface problems with error estimates independent of contrast},
	volume = {73},
	year = {2017}}

@article{2014AdjeridBenromdhaneLin,
	author = {Adjerid, Slimane and Ben-Romdhane, Mohamed and Lin, Tao},
	fjournal = {International Journal of Numerical Analysis and Modeling},
	issn = {1705-5105},
	journal = {Int. J. Numer. Anal. Model.},
	mrclass = {65N30 (65N50)},
	mrnumber = {3218337},
	number = {3},
	pages = {541--566},
	title = {Higher degree immersed finite element methods for second-order elliptic interface problems},
	volume = {11},
	year = {2014}}

@article{2023JiWangChenLi,
	author = {Ji, Haifeng and Wang, Feng and Chen, Jinru and Li, Zhilin},
	doi = {10.1051/m2an/2023047},
	fjournal = {ESAIM. Mathematical Modelling and Numerical Analysis},
	issn = {2822-7840,2804-7214},
	journal = {ESAIM Math. Model. Numer. Anal.},
	mrclass = {65N15 (35R05 65N30)},
	mrnumber = {4609882},
	number = {4},
	pages = {2041--2076},
	title = {Analysis of nonconforming {IFE} methods and a new scheme for elliptic interface problems},
	url = {https://doi.org/10.1051/m2an/2023047},
	volume = {57},
	year = {2023}}

@article{2026JiLi,
	author = {Ji, Haifeng and Li, Zhilin},
	doi = {10.1007/s10915-025-03155-x},
	fjournal = {Journal of Scientific Computing},
	issn = {0885-7474,1573-7691},
	journal = {J. Sci. Comput.},
	mrclass = {65N15 (35R05 65N30)},
	mrnumber = {5000929},
	number = {1},
	pages = {Paper No. 27, 29},
	title = {An immersed finite element method for anisotropic elliptic interface problems with nonhomogeneous jump conditions},
	url = {https://doi.org/10.1007/s10915-025-03155-x},
	volume = {106},
	year = {2026}}

@article{2022JiWangChenLi2,
	author = {Ji, Haifeng and Wang, Feng and Chen, Jinru and Li, Zhilin},
	doi = {10.1007/s00211-022-01276-1},
	fjournal = {Numerische Mathematik},
	issn = {0029-599X,0945-3245},
	journal = {Numer. Math.},
	mrclass = {65N30 (35R05 65N12 65N15)},
	mrnumber = {4405690},
	mrreviewer = {Barbara\ Verf\"urth},
	number = {4},
	pages = {1035--1086},
	title = {A new parameter free partially penalized immersed finite element and the optimal convergence analysis},
	url = {https://doi.org/10.1007/s00211-022-01276-1},
	volume = {150},
	year = {2022}}

@article{2026ZhuangZhangSarkisLin,
	author = {Zhuang, Qiao and Zhang, Zhongqiang and Sarkis, Marcus and Lin, Tao},
	doi = {10.1016/j.cam.2025.117126},
	fjournal = {Journal of Computational and Applied Mathematics},
	issn = {0377-0427,1879-1778},
	journal = {J. Comput. Appl. Math.},
	mrclass = {65N30 (65N12)},
	mrnumber = {4968689},
	pages = {Paper No. 117126, 26},
	title = {Higher-degree rectangular immersed finite elements discontinuous {G}alerkin methods for elliptic interface problems},
	url = {https://doi.org/10.1016/j.cam.2025.117126},
	volume = {476},
	year = {2026}}

@article{2019GuoLinZhuang,
	author = {Ruchi Guo and Tao Lin and Qiao Zhuang},
	journal = {International Journal Of Numerical Analysis And Modeling},
	number = {4},
	pages = {575--589},
	title = {Improved error estimation for the partially penalized immersed finite element methods for elliptic interface problems},
	volume = {16},
	year = {2019}}

@article{2025ChenZhang,
	author = {Chen, Yuan and Zhang, Xu},
	journal = {submitted},
	pages = {arXiv:2509.12555},
	title = {An immersed $C_0$ interior penalty finite element method for biharmonic interface problems.},
	year = {2025}}

@book{2008HesthavenWarburton,
	author = {Hesthaven, Jan S. and Warburton, Tim},
	doi = {10.1007/978-0-387-72067-8},
	isbn = {978-0-387-72065-4},
	mrclass = {65-02 (65M60 65N30)},
	mrnumber = {2372235},
	mrreviewer = {Weimin\ Han},
	note = {Algorithms, analysis, and applications},
	pages = {xiv+500},
	publisher = {Springer, New York},
	series = {Texts in Applied Mathematics},
	title = {Nodal discontinuous {G}alerkin methods},
	url = {https://doi.org/10.1007/978-0-387-72067-8},
	volume = {54},
	year = {2008}}

@article{2018AdjeridBenromdhaneLin,
	author = {Adjerid, Slimane and Ben-Romdhane, Mohamed and Lin, Tao},
	doi = {10.1016/j.camwa.2017.10.010},
	fjournal = {Computers \& Mathematics with Applications. An International Journal},
	issn = {0898-1221,1873-7668},
	journal = {Comput. Math. Appl.},
	mrclass = {65N30 (65N12)},
	mrnumber = {3775090},
	mrreviewer = {Rajen\ Kumar\ Sinha},
	number = {6},
	pages = {1868--1881},
	title = {Higher degree immersed finite element spaces constructed according to the actual interface},
	url = {https://doi.org/10.1016/j.camwa.2017.10.010},
	volume = {75},
	year = {2018}}

@article{2023AdjeridBabuskaGuoLin,
	author = {Adjerid, Slimane and Babu\v ska, Ivo and Guo, Ruchi and Lin, Tao},
	doi = {10.1016/j.cma.2022.115770},
	fjournal = {Computer Methods in Applied Mechanics and Engineering},
	issn = {0045-7825,1879-2138},
	journal = {Comput. Methods Appl. Mech. Engrg.},
	mrclass = {65N30 (65N12 65N15)},
	mrnumber = {4518020},
	mrreviewer = {Deepjyoti\ Goswami},
	pages = {Paper No. 115770, 37},
	title = {An enriched immersed finite element method for interface problems with nonhomogeneous jump conditions},
	url = {https://doi.org/10.1016/j.cma.2022.115770},
	volume = {404},
	year = {2023}}

@article{2022Adjerid,
	author = {Adjerid, Slimane},
	doi = {10.1016/j.camwa.2022.11.001},
	fjournal = {Computers \& Mathematics with Applications. An International Journal},
	issn = {0898-1221,1873-7668},
	journal = {Comput. Math. Appl.},
	mrclass = {65N50 (35J15 65N30)},
	mrnumber = {4512459},
	pages = {331--353},
	title = {A study of high-order immersed finite element spaces by pointwise interface conditions on curved interfaces},
	url = {https://doi.org/10.1016/j.camwa.2022.11.001},
	volume = {128},
	year = {2022}}

@article{2021ChenZhang,
	author = {Chen, Yuan and Zhang, Xu},
	journal = {Int. J. Numer. Anal. Model.},
	number = {1},
	pages = {120--141},
	title = {A {$P_2$}-{$P_1$} partially penalized immersed finite element method for Stokes interface problems.},
	volume = {18},
	year = {2021}}

@article{2024ChenZhang,
	author = {Chen, Yuan and Zhang, Xu},
	doi = {10.1007/s10915-023-02414-z},
	fjournal = {Journal of Scientific Computing},
	issn = {0885-7474,1573-7691},
	journal = {J. Sci. Comput.},
	mrclass = {65N30 (35R05 65N15)},
	mrnumber = {4675249},
	number = {1},
	pages = {Paper No. 19, 27},
	title = {Solving {N}avier-{S}tokes equations with stationary and moving interfaces on unfitted meshes},
	url = {https://doi.org/10.1007/s10915-023-02414-z},
	volume = {98},
	year = {2024}}

@article{2005KafafyLinLinWang,
	author = {Kafafy, R. and Lin, T. and Lin, Y. and Wang, J.},
	coden = {IJNMBH},
	doi = {10.1002/nme.1401},
	fjournal = {International Journal for Numerical Methods in Engineering},
	issn = {0029-5981},
	journal = {Internat. J. Numer. Methods Engrg.},
	mrclass = {78A25 (65N30 78M10)},
	mrnumber = {2172214 (2006e:78008)},
	number = {7},
	pages = {940--972},
	title = {Three-dimensional immersed finite element methods for electric field simulation in composite materials},
	url = {http://dx.doi.org/10.1002/nme.1401},
	volume = {64},
	year = {2005}}

@article{2021GuoZhang,
	author = {Guo, Ruchi and Zhang, Xu},
	fjournal = {Journal of Computational Physics},
	issn = {0021-9991},
	journal = {J. Comput. Phys.},
	mrclass = {65N30 (65N15)},
	mrnumber = {4265655},
	pages = {110445},
	title = {Solving three-dimensional interface problems with immersed finite elements: {A}-priori error analysis},
	volume = {441},
	year = {2021}}

@article{2019LinSheenZhang,
	author = {Tao Lin and Dongwoo Sheen and Xu Zhang},
	journal = {J. Sci. Comput.},
	number = {1},
	pages = {442--463},
	title = {{A nonconforming immersed finite element method for elliptic interface problems}},
	volume = {79},
	year = {2019}}

@article{2008HeLinLin,
	author = {He, Xiaoming and Lin, Tao and Lin, Yanping},
	doi = {10.1002/num.20318},
	fjournal = {Numerical Methods for Partial Differential Equations. An International Journal},
	issn = {0749-159X},
	journal = {Numer. Methods Partial Differential Equations},
	mrclass = {65N30 (35J25)},
	mrnumber = {2427191 (2009c:65307)},
	mrreviewer = {Marius Ghergu},
	number = {5},
	pages = {1265--1300},
	title = {Approximation capability of a bilinear immersed finite element space},
	url = {http://dx.doi.org/10.1002/num.20318},
	volume = {24},
	year = {2008}}

@article{2004LiLinLinRogers,
	author = {Li, Z. and Lin, T. and Lin, Y. and Rogers, R. C.},
	doi = {10.1002/num.10092},
	fjournal = {Numerical Methods for Partial Differential Equations. An International Journal},
	issn = {0749-159X},
	journal = {Numer. Methods Partial Differential Equations},
	mrclass = {65N30 (41A25 65N15)},
	mrnumber = {2046521 (2005f:65153)},
	mrreviewer = {Josef Dal{\'{\i}}k},
	number = {3},
	pages = {338--367},
	title = {An immersed finite element space and its approximation capability},
	url = {http://dx.doi.org/10.1002/num.10092},
	volume = {20},
	year = {2004}}

@article{2015LinLinZhang,
	author = {Lin, Tao and Lin, Yanping and Zhang, Xu},
	doi = {10.1137/130912700},
	fjournal = {SIAM Journal on Numerical Analysis},
	issn = {0036-1429},
	journal = {SIAM J. Numer. Anal.},
	mrclass = {65N30 (35R05 65N15)},
	mrnumber = {3338673},
	number = {2},
	pages = {1121--1144},
	title = {Partially penalized immersed finite element methods for elliptic interface problems},
	url = {http://dx.doi.org/10.1137/130912700},
	volume = {53},
	year = {2015}}

@article{2017AdjeridGuoLin,
	author = {Adjerid, Slimane and Guo, Ruchi and Lin, Tao},
	fjournal = {International Journal of Numerical Analysis and Modeling},
	issn = {1705-5105},
	journal = {Int. J. Numer. Anal. Model.},
	number = {4-5},
	pages = {604--626},
	title = {High Degree Immersed Finite Element Spaces by A Least Squares Method},
	volume = {14},
	year = {2017}}

@article{2019GuoLin2,
	author = {Guo, Ruchi and Lin, Tao},
	doi = {10.1137/18M121318X},
	fjournal = {SIAM Journal on Numerical Analysis},
	issn = {0036-1429},
	journal = {SIAM J. Numer. Anal.},
	mrclass = {65N30 (35J25 35R05 65N12)},
	mrnumber = {3977109},
	number = {4},
	pages = {1545--1573},
	title = {A higher degree immersed finite element method based on a {C}auchy extension for elliptic interface problems},
	url = {https://doi.org/10.1137/18M121318X},
	volume = {57},
	year = {2019}}

@article{2025AdjeridLinMeghaichi,
	author = {Adjerid, Slimane and Lin, Tao and Meghaichi, Haroun},
	doi = {10.1016/j.cma.2025.117829},
	fjournal = {Computer Methods in Applied Mechanics and Engineering},
	issn = {0045-7825,1879-2138},
	journal = {Comput. Methods Appl. Mech. Engrg.},
	mrclass = {65N30 (65N12)},
	mrnumber = {4867707},
	pages = {Paper No. 117829, 19},
	title = {The {F}renet immersed finite element method for elliptic interface problems: an error analysis},
	url = {https://doi.org/10.1016/j.cma.2025.117829},
	volume = {438},
	year = {2025}}

@article{2024AdjeridLinMeghaichi,
	author = {Adjerid, Slimane and Lin, Tao and Meghaichi, Haroun},
	doi = {10.1016/j.cma.2023.116703},
	fjournal = {Computer Methods in Applied Mechanics and Engineering},
	issn = {0045-7825,1879-2138},
	journal = {Comput. Methods Appl. Mech. Engrg.},
	mrclass = {65N30 (65N12)},
	mrnumber = {4682175},
	pages = {Paper No. 116703, 21},
	title = {A high order geometry conforming immersed finite element for elliptic interface problems},
	url = {https://doi.org/10.1016/j.cma.2023.116703},
	volume = {420},
	year = {2024}}

@article{2025AdjeridLinMeghaichi2,
	author = {Adjerid, Slimane and Lin, Tao and Meghaichi, Haroun},
	journal = {arXiv preprint},
	pages = {arXiv:2510.12018},
	title = {Construction of Basis Functions for the Geometry Conforming Immersed Finite Element Method},
	year = {2025}}

@article{2019GuoLin,
	author = {Ruchi Guo and Tao Lin},
	journal = {IMA J. Numer. Anal.},
	number = {1},
	pages = {482--511},
	title = {{A group of immersed finite element spaces for elliptic interface problems}},
	volume = {39},
	year = {2019}}

@article{1998Li,
	author = {Li, Zhilin},
	coden = {ANMAEL},
	doi = {10.1016/S0168-9274(98)00015-4},
	fjournal = {Applied Numerical Mathematics. An IMACS Journal},
	issn = {0168-9274},
	journal = {Appl. Numer. Math.},
	mrclass = {65L20 (65N15)},
	mrnumber = {1634348 (99h:65140)},
	mrreviewer = {Michael Fr{\"o}hner},
	number = {3},
	pages = {253--267},
	title = {The immersed interface method using a finite element formulation},
	url = {http://dx.doi.org/10.1016/S0168-9274(98)00015-4},
	volume = {27},
	year = {1998}}

@book{2012AbateTovena,
	author = {Abate, Marco and Tovena, Francesca},
	publisher = {Springer Science \& Business Media},
	title = {Curves and surfaces},
	year = {2012}}

%
%
%
%
%
%
%
%
\end{document}